\documentclass[10pt]{article}

\usepackage{mathrsfs}
\usepackage{amsfonts}
\usepackage{}

\usepackage[leqno]{amsmath}
\usepackage{graphicx}
\usepackage{latexsym}
\usepackage{amsmath,amsfonts,amssymb,amsthm,mathrsfs,euscript,makeidx,color}
\usepackage{enumerate}
\def\BBsi{\boldsymbol\sigma}

\allowdisplaybreaks[2]

\def\sqr#1#2{{\vcenter{\vbox{\hrule height.#2pt
				\hbox{\vrule width.#2pt height#1pt \kern#1pt \vrule width.#2pt}
				\hrule height.#2pt}}}}
\def\5n{\negthinspace \negthinspace \negthinspace \negthinspace \negthinspace }
\def\4n{\negthinspace \negthinspace \negthinspace \negthinspace }
\def\3n{\negthinspace \negthinspace \negthinspace }
\def\2n{\negthinspace \negthinspace }
\def\1n{\negthinspace }

\def\dbE{\mathbb{E}}
\def\dbF{\mathbb{F}}

\def\dbH{\mathbb{H}}

\def\dbM{\mathbb{M}}

\def\dbP{\mathbb{P}}

\def\dbR{\mathbb{R}}

\def\sC{\mathscr{C}}
\def\sD{\mathscr{D}}

\def\sL{\mathscr{L}}

\def\sN{\mathscr{N}}

\def\sP{\mathscr{P}}

\def\sU{\mathscr{U}}

\def\cA{{\cal A}}

\def\cF{{\cal F}}

\def\cH{{\cal H}}

\def\cM{{\cal M}}
\def\cN{{\cal N}}

\def\cP{{\cal P}}

\def\cT{{\cal T}}
\def\cU{{\cal U}}

\def\BJ{{\bf J}}

\def\BV{{\bf V}}

\def\Bb{{\bf b}}

\def\Bg{{\bf g}}
\def\Bh{{\bf h}}

\def\Bw{{\bf w}}

\def\ds{\displaystyle}

\def\ns{\noalign{\ss}}

\def\q{\quad}
\def\qq{\qquad}

\def\({\Big (}
\def\){\Big )}
\def\[{\Big[}
\def\]{\Big]}

\def\lan{\langle}
\def\ran{\rangle}

\def\a{\alpha}

\def\g{\gamma}
\def\d{\delta}
\def\e{\varepsilon}

\def\l{\lambda}

\def\si{\sigma}
\def\t{\tau}

\def\th{\theta}

\def\Th{\Theta}
\def\L{\Lambda}

\def\wt{\widetilde}

\def\cd{\cdot}

\def\bde{\begin{definition}\label}
	\def\ede{\end{definition}}
	\def\bel{\begin{equation}\label}
		\def\ee{\end{equation}}
	\def\bt{\begin{theorem}\label}
		\def\et{\end{theorem}}
	\def\bc{\begin{corollary}\label}
		\def\ec{\end{corollary}}
	\def\bl{\begin{lemma}\label}
		\def\el{\end{lemma}}
	\def\bp{\begin{proposition}\label}
		\def\ep{\end{proposition}}
	\def\bas{\begin{assumption}\label}
		\def\eas{\end{assumption}}
	\def\br{\begin{remark}\label}
		\def\er{\end{remark}}
	\def\bex{\begin{example}\label}
		\def\ex{\end{example}}
	\def\ba{\begin{array}}
		\def\ea{\end{array}}
	\def\ben{\begin{enumerate}}
		\def\een{\end{enumerate}}
	
	\def\square#1{\vbox{\hrule\hbox{\vrule height#1%
				\kern#1\vrule}\hrule}}
	\def\rectangle#1#2{\vbox{\hrule\hbox{\vrule height#1%
				\kern#2\vrule}\hrule}}

	\font\tenbb=msbm10 \font\sevenbb=msbm7 \font\fivebb=msbm5
	
	\newfam\bbfam
	\scriptscriptfont\bbfam=\fivebb \textfont\bbfam=\tenbb
	\scriptfont\bbfam=\sevenbb

	\newtheorem{theorem}{\indent Theorem}[section]
	\newtheorem{definition}[theorem]{\indent Definition}
	\newtheorem{proposition}[theorem]{\indent Proposition}
	\newtheorem{corollary}[theorem]{\indent Corollary}
	\newtheorem{lemma}[theorem]{\indent Lemma}
	\newtheorem{remark}[theorem]{\indent Remark}
	\newtheorem{example}[theorem]{\indent Example}
	
	\newtheorem{assumption}[theorem]{\indent Assumption}
	
	\makeatletter
	
	\@addtoreset{equation}{section}
	\makeatother
	
	\def\bea{\begin{equation*}}
		\def\eea{\end{equation*}}
	
	\def\bel{\begin{equation}\label}
		\def\eel{\end{equation}}
	
	\def\ba{\begin{array}}
		\def\ea{\end{array}}
	\newcommand{\ad}{&\!\!\!\displaystyle}

	\def\({\Big (}
	\def\){\Big )}
	\def\[{\Big[}
	\def\]{\Big]}
	\def\q{\quad}
	\def\qq{\qquad}
	
	\def\d{\delta}
	\def\e{\varepsilon}
	
	\def\ds{\displaystyle}

	\def\ns{\noalign{\smallskip}}

	\def\argmin{\mathop{\rm argmin}}
	\def\law{\mathop{\rm law}}

\begin{document}

\title{Closed-loop Equilibria for	Mean-Field Games
	in Randomly Switching Environments with General Discounting Costs\thanks{H. Mei was partially supported by Simons Foundation's Travel Support for Mathematicians Program (No. 00002835); G. Yin  was supported in part by the National Science Foundation under grants		DMS-2204240.}}

\author{Hongwei Mei,\thanks{Department of Mathematics and Statistics, Texas Tech University, Lubbock, TX 79409. Email: hongwei.mei@ttu.edu}~~~~
	 Son Luu Nguyen\thanks{Department of Mathematical Sciences, Florida Institute of Technology, Melbourne, FL 32901
		Email: snguyen@fit.edu}
	~~~~and~~~~ George Yin\thanks{Department of Mathematics, University of Connecticut, Storrs, CT 06269. Email: gyin@uconn.edu}}

	\maketitle
	
	\begin{abstract}
		This work is devoted to finding the closed-loop equilibria for a  class of mean-field games (MFGs)  with infinite many symmetric players in a common switching environment  when the cost functional is under
general discount in time. There are two key challenges in the application of the well-known Hamilton-Jacobi-Bellman and Fokker-Planck (HJB-FP)
approach to our problems: the path-dependence  due to the conditional mean-field interaction and the time-inconsistency  due to the general discounting cost.
To overcome the  difficulties,
a theory for a class of systems of path-dependent equilibrium  Hamilton-Jacobi-Bellman equations (HJBs) is developed. Then 
closed-loop equilibrium strategies can be identified through  a  two-step verification procedure. It should be noted that
the  closed-loop equilibrium strategies obtained satisfy a new form of local optimality in the Nash sense.
The theory obtained extends the HJB-FP approach for classical MFGs to
more general conditional MFGs  with general discounting costs.
	\end{abstract}
	

	\noindent{\bf MSCcodes.}
		60H10, 91A16, 93E03,  34K50.

	\section{Introduction}
	Since the publications of the seminal works of
	Huang, Caines, and Malham\'e \cite{Huang2003, Huang2007},
and Lasry and Lions \cite{Lasry2006a,Lasry2006b}
(published independently), the mean-field game (MFG) theory
has drawn increasing attention in the last two decades.
In a mean-field system, 	each agent (also known as a player)  plays an insignificant role, however the system as a whole is significantly impacted by the agents' combined efforts.
Bensoussan, Frehse, and Yam
\cite{BFY13} provided an illuminating
discussion of certain aspects of mean-field games and  mean-field type controls,
describing their similarities and differences together with a unified approach for treating them. For a
survey, see also \cite{Gomes} and references therein. For the most recent developments, we refer to \cite{Barreiro2019,Bayraktar2018a,Bensoussan2017,
Bensoussan2020,
Nutz2018,Sen2019} among others.
A comprehensive treatment of mean-field game theory can be found in the recent book of Carmona and  Delarue \cite{CD18}.

Along another line,  hybrid systems have gained increasing popularity owing to their ability to handle numerous real-world applications in which discrete events and continuous dynamics coexist and interact. One of the commonly used models is
switching diffusion.
The switching process  is used to depict
random jumps
 that can be modeled using a continuous-time Markov chain;
see \cite{Yin2009} for switching diffusions and applications.
For a wide range of applications, we mention the work on controlled piecewise deterministic Markov
processes \cite{Costa2010}, population dynamics,
mathematical biology
and ecology \cite{NNY21, NNY21a},
among others.
Furthermore, there were also
efforts in treating switching diffusions in conjunction with mean-field interactions \cite{Xi2009}.

The current paper aims to  find closed-loop equilibria for
a class of MFGs with infinitely-many
symmetric players sharing  a common switching environment.
Different from classical MFGs, the mean-field interaction takes a form of conditional expectation in this paper (see \cite{NguyenYinHoang 2020}) and our problems belong to a wider class of the conditional MFGs.

{In general, the equilibria for (conditional) MFGs can be classified into two  categories, namely open-loop and closed-loop.  The open-loop equilibria take the form of stochastic processes that depend on the initial time and the initial state of an agent (see, for example, \cite {JLSY24,Motte2022}). In contrast,  closed-loop equilibria are deterministic strategies represented as feedback functions of the observations. Our paper focuses on the latter one.
Moreover,
we will consider  more complex MFGs with generalized  discounting cost functional in time (compared to the special exponential discounting). The different discounting preference originated
from different subjectivity of people’s preference on the future risks in the decision process leads to a
wider range of applications in different areas (see \cite{Ekeland2010, Yong2012} for more details).}

To find the closed-equilibria, there are two key challenges in the application of the well-known  HJB-FP
(Hamilton-Jacobi-Bellman and Fokker-Planck)
approach to our problems: the path-dependence structure due to the conditional mean-field interaction and the time-inconsistency  due to the general discounting cost.
	To overcome the first challenge, a  path-dependence theory for HJBs is needed.  Due to the second feature, the optimal control problem in the HJB-step becomes time-inconsistent. To overcome this challenge,  we  introduce  the equilibrium HJB inspired by  J. Yong (see \cite{Yong2012}) to  seek a time-consistent closed-loop equilibrium  which  satisfies some local optimality only.   Combining those two parts,  our main efforts  are devoted to developing  a theory for a class of systems of  path-dependent equilibrium HJBs such that a new type of local-optimal closed-loop equilibrium strategies can be identified through  a  two-step verification procedure.
	The cases without switching were treated in  \cite{Mei2020}, whereas the current paper focuses on hybrid systems and the path-dependence yields further technical challenges.
Let us highlight the main contributions of the current paper.
\begin{itemize}
\item[(a)] Using a conditional mean field, we formulate the controlled switching diffusion with mean-field interactions.
\item[(b)] In addition to the mean-field interactions, we also take into account general discounting costs, which result in time inconsistency. The time-inconsistency leads to a new local version of Nash's optimality.
\item[(c)] In light of (a), our primary focus is on resolving the issues brought on by path dependence.
We develop novel methods for solving  a class of
systems of
 path-dependent equilibrium HJBs on the path space of the Markov chain that is
piecewise constant. In contrast to the use of
viscosity solutions
in the   previous works (e.g., \cite{EKren2016a,EKren2016b}),
because of
the special structure, we are able to find a classical solution in an appropriate sense by the fundamental solution method in partial differential equations (PDEs) \cite{Friedman}. Some critical and technical
estimates are derived
as a by-product, which are interesting in their own right.
\item[(d)] Based on the theory of path-dependent HJB equations
that we develop,
we prove the existence and uniqueness of a closed-loop equilibrium strategy for the 
MFG
sharing the same  environment and with a general discounting cost.
\end{itemize}

\begin{example}\label{m-exm}{\rm
	Before proceeding further, let us begin with a motivational
	heterogeneous agent model  in macroeconomics; see \cite{Huagge1993}.  Suppose that there are
	$N$ agents,  that
	$\a$ is a finite-state Markov chain with state space $M=\{1,\dots, m\}$ representing the aggregate shocks from the market, that
	$w_t^i: i=1,2\dots$  are independent Brownian motions being independent of  $\a$, which
	represent the heterogeneous market fluctuations for different agents,
	that  the wealth  of agent $i$
	is denoted by $X_t^i$, and   the   average wealth process is defined  by
	$
	\Pi^N_t=\frac1N\sum_{i=1}^N X_t^{i,N}.
	$
	Suppose the wealth of the agent $i$ satisfies
	\begin{equation}\label{SDEEx}
		dX_t^{i,N}= b(t,X_t^{i,N},\a_t, \Pi^{N}_t;u^i_t)dt+\si(t,X_t^{i,N})d w_t^i.
	\end{equation}
	The agent $i$ takes the action $u^i:=\{u_t^i:t\in[0,T]\}$  to maximize a profit functional as follows
	$$J^i(t,x;u^i)=\dbE\(\int_t^T e^{-\lambda (s-t)}g(X_s^{i,N},\a_s, \Pi^{N}_s;u^i_s)ds+e^{-\lambda (T-t)}h(X_T^{i,N},\a_T, \Pi^{N}_T)\). $$
	The interest to  find the appropriate behavior of each agent
	when the system consists of a large number of agents.
	Letting $N\rightarrow \infty$,  by the law of large numbers,  it can be proved that $\Pi_t^N\to \dbE [X_t|\cF_{t^-}^\a]$ as $N\to \infty$, where $X_t$ is a ``representative'' of $X_t^i$ and
	$\cF_{t^-}^\a$ is the  filtration of $\a$ up to time $t^-$. Note that all the agents are symmetric because $b,\sigma,f,g$ are the same for different agents. The superscript $i$ can be omitted in the sequel.
	Denote by
	$X^{\xi,u}$  the solution to  \eqref{SDEEx1} with initial $X_0=\xi$ (the wealth distribution of all agents initially) and
	the average wealth process by $\eta^{\xi,u}_s:=\dbE [X^{\xi,u}_t|\cF_{t^-}^\a]$.
	Then,  the heterogeneous agent model with infinite-many symmetric players reduces to optimizing
	\begin{equation}\label{costEx1}J(t,x;\eta^{\xi,u},u)=\dbE\(\int_t^Te^{-\lambda (s-t)}g(s,X_s,\a_s, \eta^{\xi,u}_s;u_s)ds+e^{-\lambda (T-t)}h(X_T,\a_T, \eta^{\xi,u}_T)\),\end{equation}
	subject to
	\begin{equation}\label{SDEEx1}
		dX_t= b(t,X_t,\a_t, \eta_t^{\xi,u};u_t)dt+\si(t,X_t)d w_t.
	\end{equation}
	Such a problem (especially when $\a $ is a constant process) has been well studied in  the mean-field game theory where the main object is to find an equilibrium  strategy $u^*$  (under suitable conditions) in Nash's sense  such that
	\bel{optmaNash}J(t,x;\eta^{\xi,u^*},u)\leq J(t,x;\eta^{\xi,u^*},u^*)\text{ for any strategy }u.\eel
	Here $\eta^{\xi,u^*}$ is considered as a pre-committed average wealth process if all the agents take the same ``rational"  strategy $u^*$. The optimality in \eqref{optmaNash} essentially says that the equilibrium strategy is optimal if all  other agents adopt such a strategy. We emphasize that our problem is different from a McKean-Vlasov optimal control problem where an optimal control $u^*$ satisfies
 $J(t,x;\eta^{\xi,u},u)\leq J(t,x;\eta^{\xi,u^*},u^*)$ for any strategy  $u$.
		The exponential discounting structure in the profit functional in \eqref{costEx1}  essentially leads to a time-consistent optimal control problem for each agent to solve. In reality,
	it is often necessary to consider
	some general discounting profits functionals.  That is, the discounting factors might be  replaced by some other functions $\mu(s, t)$ and $\nu(T, t)$.
	The so-called hyperbolic discounting
is such an example
 by letting $\mu(s, t)=\frac1{1+\lambda(s-t)},\q \nu(T, t)=\frac1{1+\lambda(T-t)}.$
	For general $\mu$ and $\nu$, the optimization  problem for each agent is time-inconsistent;
  the dynamic programming principle fails.
	For more details, the reader is referred to  Section 2 in \cite{Yong2012}. In this paper, we will consider the following profit-functional
	$V(t,x;\eta^{\xi,u},u)=\wt J(t;t,x;\eta^{\xi,u},u)$ where
	\begin{equation}\label{costEx2}\wt J(\tau;t,x;\eta^{\xi,u},u)=\dbE\(\int_t^T\wt g(\t;s,X_s,\a_s, \eta^{\xi,u}_s;u_s)ds+\wt h(\t,X_T,\a_T, \eta^{\xi,u}_T)\).\end{equation}
	Because of the same reason, it is impossible to find an equilibrium strategy satisfying \eqref{optmaNash}.  {Instead, a local version of   \eqref{optmaNash} will be proved in the future
to illustrate the local optimality of the  new equilibrium strategy in our paper (see \eqref{localopti}).}
	}\end{example}

	The rest of the paper is arranged as follows.  Section \ref{sec:for} presents
	the formulation of our problem. In addition we introduce the path space of the Markov chain and work out the
	calculus on the path space. To go through the HJB-FP
method in our problems, we  detail
the  FP-step and HJB-step in Section \ref{sec:fp} and Section \ref{sec: HJB-step}, respectively.
	Especially, a  theory on  path-dependent
systems of
 equilibrium HJBs is developed, which presents  the main contribution of the paper.
	Using the results of Section \ref{sec: HJB-step},  we proceed with our final goal--to obtain the closed-loop equilibria for  MFGs with switching under general discounting in Section \ref{sec: mfg}.
	We also revisit  the motivational example on  heterogeneous agent model in   Section \ref{sec: mfg} to illustrate the theory developed.  Finally, concluding remarks are made in Section \ref{sec:cr}.
	
		\section{Formulation}\label{sec:for}
	On a complete probability space $(\Omega,\dbF,\dbP)$ with $\cN$ being the set of all $\dbP$-null sets, let  $\alpha(\cdot)$ be a
	finite-state Markov chain with state space  $M=\{1,\cdots,m\}$ and  generator $Q=(q_{ij})_{m\times m}$. Denote by $\cF_{t^-}^\alpha$ the smallest $\sigma$-algebra containing $\{\alpha(s):0\leq s< t\}$ augmented with all $\dbP$-null sets, and $\dbF$ be  the natural filtration of $\a(\cdot)$ augmented with all $\dbP$-null sets. Assume that $\cF_{0^-}^\a=\{\emptyset, \Omega\}\bigvee \cN=\cF_{0}^\a$. Let $B(\cdot)$ be a standard $n$-dimensional Brownian motion (independent of $\a$) and $\cF^B$ be its natural filtration  augmented with all $\dbP$-null sets.  Denote $\cF_t:=\cF_{t}^\alpha \bigvee \cF_t^B\bigvee \cN$, $\cF_{t^-}:=\cF_{t^-}^\alpha \bigvee \cF_t^B\bigvee \cN$. 	Let $\sP(\dbR^n)$ ($\sP_2(\dbR^n)$, resp.) be the space of probability measures (with finite second moment, resp.) on $\dbR^n$. It is assumed that   there are infinite many symmetric players.
We  formulate the conditional MFGs as  controlled   McKean-Vlasov dynamics with switching.

 Consider
	\bel{SDE-1}\left\{\ba{ll}\ns\ad dX(t)=\Bb(t,\alpha(t^-),X(t),\eta(t);u(t))dt+\BBsi(t, X(t))dB(t),\\
	\ns\ad\eta(t)=\law(X(t)|\cF_{t^-}^\alpha),~ X(0)=\xi,~\a(0^-)=i,\ea\right.\eel
	where $\Bb:[0,T]\times M\times \dbR^n\times \sP(\dbR^n)\times U\mapsto \dbR^n$, $\BBsi:[0,T]\times \dbR^n\mapsto \dbR^{n\times n}$ are appropriate maps.  The $B(\cdot)$ is a standard $n$-dimensional Brownian motion.
	Especially
	$\eta(t)$   is the  conditional probability law of  $X(t)$ on the filtration $\cF_{t^-}^\alpha$ in the sense that for any bounded and  continuous function $f:\dbR^n\mapsto \dbR$,
	$\dbE[f(X(t))|\cF_{t^-}^\alpha]=\int_{\dbR^n}f(x)\eta(t, dx).$  In our problem, the player aims to find an appropriate $u(\cdot)\in \cU[0,T]$ in accordance with the cost functional   \bel{cost}\BV(t,\xi,i;u(\cdot))=\BJ(t;t,\xi,i;u(\cdot)).\eel
 where
\begin{align*}&\BJ(\t;t,\xi,i;u(\cdot))\\
&\q=\dbE_{t,
	\xi,i}\(\int_t^T \Bg(\t;s,\alpha(s^-),X(s),\eta(s);u(s))ds +\Bh(\t;\alpha(T^-),X(T),\eta(T))\), \\
&\cU[0,T]:=\Big\{u:[0,T]\times\Omega\mapsto U\big| u(\cdot) \text{ is $\dbF$-adapted with }\dbE\int_0^T|u(s)|^2ds<\infty\Big\}.\nonumber\end{align*}
 The set $U$ is the action space, and 	 $\Bg:[0,T]\times[0,T]\times M\times \dbR^n\times \sP(\dbR^n)\times U\mapsto \dbR^+$ and $\Bh:[0,T]\times M\times\dbR^n\times \sP(\dbR^n)\mapsto\dbR^+$ are two appropriate functions denoting the running cost and terminal cost respectively. Note that the additional time factor $\tau$ yields a general discounting  cost in our problem. As was  mentioned,
we aim to find  closed-loop equilibria  instead of  optimal controls for
the  McKean-Vlasov dynamics. We will follow the well-known HJB-FP approach while some substantial modifications are needed for  our problems. The following is the main assumption for our problem. Let $\sP_2(\dbR^n)\subset \sP(\dbR^n)$ be equipped with the Wasserstein 2 metric $\Bw$ (see \cite{Villani 2008}).
\begin{assumption}\label{A-3}
\rm (1) Suppose there exist $b_1:[0,T]\times  M\times \dbR^n\times  U\mapsto\dbR^n, b_2:[0,T]\times  M\times \dbR^n\times \cP(\dbR^n)\mapsto\dbR^n$ and $g_1:[0,T]\times[0,T]\times  M\times \dbR^n\times  U\mapsto\dbR^n, g_2:[0,T]\times[0,T]\times  M\times \dbR^n\times \cP(\dbR^n)\mapsto\dbR^n$ such that $\Bb(t,i,x,\rho;v)=b_1(t,i,x;v)+b_2(t,i,x,\rho)$, $\Bg(\t;t,i,x,\rho;v)=g_1(t,i,x;v)+g_2(\tau;t,i,x,\rho).$

(2) The action space $U$ is a compact subset of $\dbR$ and  there exists a mapping $\psi:[0,T]\times  M\times  \dbR^n\times  \dbR^n\mapsto U$ defined by
\begin{equation}\label{defpsi}\psi(t,i,x,p)=\argmin_{v\in U}\big( \langle p, b_1(t,i,x;v)\rangle+g_1(t,i,x;v)\big), \end{equation}
which is continuous w.r.t. $t$ and uniformly Lipschitz w.r.t. $x$ and $p$  with Lipschitz constant $L_\psi$ and  $|\psi(t,i,0,v)|<K$.

(3)
 $b_1$ and $g_1$ are continuous w.r.t. $t$, and $g_2, \Bh$ are continuous w.r.t. $(\tau,t)$ satisfying  $$\left\{\ba{ll}\ns\ad |b_1(t,i,x;v)|+|b_2(t,i,x,\rho)|+|g_1(t,i,0;0)|+|g_2(\t;t,i,0,\rho)|\leq K,\\
\ns\ad |b_1(t,i,x_1;v_1)-b_1(t,i,x_2;v_2)|+|g_1(t,i,x_1;v_1)-g_1(t,i,x_2;v_2)|\\
\ns\ad\qq\qq\leq K(|x_1-x_2|+|v_1-v_2|),\\

\ns\ad |b_2(t,i,x_1,\rho_1)-b_2(t,i,x_2,\rho_2)|+|g_2(\tau;t,i,x_1,\rho_1)-g_2(\t;t,i,x_2,\rho_2)|\\
\ns\ad\q+|\Bh(\t;i,x_1,\rho_1)-\Bh(\t;i,x_2,\rho_2)\leq K(|x_1-x_2|+\Bw(\rho_1,\rho_2)),\\
\ns\ad |\Bh(\t,t, 0,\rho)|+|D_{xx}\Bh|+|\BBsi(s,0)|+|D_x\BBsi(s,x)|\leq K,\\
\ns\ad \l_0^{-1}I\leq \BBsi(s,x)\BBsi^\top(s,x)\leq \l_0I\text{ for some }\l_0>1.\ea \right.$$

\end{assumption}

 From Assumption \ref{A-3}, we see that  $v$ and $\rho$ are separated in $\Bb$ and $\Bg$, which is a general assumption for mean-field game theory (see \cite{Ca2015} for example). In this scenario, we have a separated Hamiltonian $\cH$ defined by
\begin{align*}&\cH(\tau; t,i,x,\rho, p,q):=\langle p, b_1(t,i,x;\psi(t,x,i,q))\rangle+g_1(t,i,x;\psi(t,x,i,q))\\
&\qq+\langle p, b_2(t,i,x,\rho)\rangle+g_2(\tau;t,i,x,\rho).
\end{align*} In fact, if $\Bb$ and $\Bh$ are independent of $\tau$, the Hamiltonian can be simplified by letting $p=q$, which reduces to the classical case
\begin{align*}&\bar\cH( t,i,x,\rho, p)=\min_{v\in U}\(\langle p, b_1(t,i,x;v)\rangle+g_1(t,i,x;v)\)+\langle p, b_2(t,i,x,\rho)\rangle+g_2(t,i,x,\rho).
\end{align*}
The  introduction of $\tau$ and $q$ in $\cH$ is inspired by \cite{Yong2012} to tackle the time-inconsistency. With all those notations, our modified HJB-FP approach can be summarized as follows.

(1) {\it FP-step}: Given a pre-committed feedback strategy $u$ depending on time, state, and the path of the Markov chain, the solution to  \eqref{SDE-1} leads to a  map $\zeta(t,\a_t)$ representing   the conditional distribution of $X(t)$ on $\cF_{t^-}^\a$.

(2) {\it HJB-step}: Replace $\eta(t)$ by the feedback form of $\zeta(t,\a_t)$ where $\a_t$ is the path of $\a(\cd)$ up to time $t^-$. Then we  solve a  control problem subject to \eqref{SDE-1} under the general discounting cost \eqref{cost}. Different from the classical MFGs,  our control problem turns out to be time-inconsistent and  path-dependent. To work with the two features,  
we develop a theory for systems of path-dependent equilibria of HJBs  \begin{align}\label{PHJE-10000}\begin{cases}&\!\!\!\!\!(\partial^\alpha \!+\!\cA)\Theta(\t;t,w_t,x)\!\!+\!\!\cH(\tau;t,\! w_t(t^-),\! x,\!\zeta(t,w_t),\!D_x\Theta(\t;t,w_t,x),\!D_x\Theta(t;t,w_t,x))\!=\!0,\\
& \!\!\!\!\!\Theta(\t;T,w_T,x)=\Bh(\t;w_T(T^-),x,\zeta(t,w_T));\end{cases}\end{align}
where  the operator $\cA\Th(\t;t,w_t,x):=\text{Trace}[a(t,x) D_x^2\Th(\t;t,w_t,x)]/2$ and  $\partial^\a$ is some   appropriate path-derivative to be defined later. The   $w_t$ denotes a path of $\a$ upto time $t^-$ and $w_t(s)=\a(s)$ for $s< t$.
If \eqref{PHJE-10000} can be solved appropriately, one can   derive  a closed-loop control $u^*(t,\a_t,X(t))=\psi(t, \a(t^-),X(t), D_x\Th(t,\a_t,X_t))$. 

(3) {\it Closed-loop equilibrium}: If $u^*$ coincides with the pre-committed $u$, the feedback strategy $u$ is the required closed-loop equilibrium whose exact definition will be given later in Section \ref{sec: mfg}.

From the discussion above,  a closed-loop equilibrium is essentially a fixed point in a two-step verification procedure, where the key challenge  lies in developing  a theory on the path-dependent
 HJBs \eqref{PHJE-10000}. This  presents the main novelty of the paper.
Moreover,   we will derive  a  new local optimality (comparable to \eqref{optmaNash}) for the closed-loop equilibrium obtained. In particular, if the cost functional in \eqref{cost} is exponential discounting,  the local optimality is equivalent to   \eqref{optmaNash}. This further signifies that our theory is a generalization of the classical HJB-FP approach for MFGs with switching under general discounting costs.
Now, it is natural to start our investigation on the path space of $\a$ and the path-derivative $\partial^\a$. More details will be presented in the following subsection.



\subsection{Path Space and the Calculus}
In this subsection, we develop the calculus on  the path space of  Markov chain $\alpha$. Write  $\alpha_t:=\{\alpha(s):0\leq s< t\}$ by the path of $\a(\cdot)$ up to  time $t^-$. Note that all the paths of $\a(\cdot)$ are  piecewise constant and c\'adl\'ag (right continuous with left limit).
We equip $M$, the state space of $\a$,  with  the discrete metric
$d_M(i,j)=\delta_{ij}$, where $\delta_{ij}$ is the Kronecker delta function.

Write $w_{[t_1,t_2)}:[t_1,t_2)\mapsto M $ by
a possible path of $\a(\cd)$ on time interval $[t_1,t_2)$.
The set of all possible paths of $\a(\cd)$ on $[t_1,t_2)$ is defined by
$$\ba{ll}\cM[t_1,t_2):=\{w_{[t_1,t_2)}(\cdot):[t_1,t_2)\mapsto M|~w_{[t_1,t_2)}(\cdot) \text{ is c\'adl\'{a}g}\}.\ea$$
For any $w_{[t_1,t_2)}\in \cM[t_1,t_2)$, we can continuously extend the path to the right-endpoint $t_2$ by $$w^+_{[t_1,t_2]}(s):=\left\{\ba{ll}\ad w_{[t_1,t_2)}(s)\qq\text{ for }s\in[t_1,t_2);\\
\ns\ad w_{[t_1,t_2)}(t_2^-)\q\text{ for }s=t_2.\ea\right.$$
Such an observation allows us to define  the Skorohod-type metric $\sD$ for  $\cM[t_1,t_2)$   by
$$\sD(w_{[t_1,t_2)},\tilde w_{[t_1,t_2)})=\inf_{\lambda\in\L}\sup_{t_1\leq s\leq t_2}\(|\l(s)-s|+d_M(w^+_{[t_1,t_2]}(\l(s)),\tilde w^+_{[t_1,t_2]}(s))\),
\ \hbox{ where}$$
$$
\ba{ll}\ad \L:=\Big\{\l:[t_1,t_2]\mapsto[t_1,t_2]\Big|~\l\text{  is continuous and strictly increasing}\\
\ad\qq\qq\qq \text{ with $\l(t_1)=t_1$ and $\l(t_2)=t_2$}\Big\}.\ea
$$
We denote by
$\sN(w_{[t_1,t_2)})$ the number of jumps in $w_{[t_1,t_2)}.$
Note that the definitions of $\sD$ and $\sN$ depend  on the time-interval $[t_1,t_2)$, while we will omit such dependence in the notations for the sake of simplicity. The  combination of  two paths $w_{[t_0,t_1)}$ and $w_{[t_1,t_2)}$ for $t_0<t_1<t_2$  is defined by  $$[w_{[t_0,t_1)}\oplus w_{[t_1,t_2)}] (t):=\left\{\ba{ll}w_{[t_0,t_1)}(t)\text{ for }t\in [t_0,t_1),\\ [2mm]
w_{[t_1,t_2)}(t)\text{ for }t\in [t_1,t_2).\ea\right.$$
For simplicity, we write $\cM_t= \cM[0,t)$, $\cM=\cM[0,T) $,
and $w_{[0,t)}=w_t.$ Moreover, for $w\in\cM$, write $w_t$
as the restriction of $w$ on $[0,t)$.
We have the following observations for $\sD$ of which the proof is obvious.

\begin{proposition}
{\rm (1)}
If  $\sD(w_{[t_1,t_2)},\tilde w_{[t_1,t_2)})<1$,
then
$\sN(w_{[t_1,t_2)}=\sN(\tilde w_{[t_1,t_2)})$, and
\bel{SNSD}\int_{t_1}^{t_2}I(w_{[t_1,t_2)}(t)\neq \tilde w_{[t_1,t_2)}(t))dt\leq \sN(w_{[t_1,t_2)})\sD(w_{[t_1,t_2)},\tilde w_{[t_1,t_2)}).\eel

{\rm (2)} It follows that $$\sD(\g_{[t_0,t_1)}\oplus w_{[t_1,t_2)},\g_{[t_0,t_1)}\oplus\tilde w_{[t_1,t_2)})\leq  \sD( w_{[t_1,t_2)},\tilde  w_{[t_1,t_2)}).$$

\end{proposition}

Define time-augmented path spaces of $\alpha(\cd)$ as follow
\begin{align*}&\dbM:=\bigcup_{0\leq t\leq T}\{t\}\times \cM_t,\\
& \dbM_{\gamma_t}:=\{(s, w_s)\in\dbM: s\in[t,T]\text{ and }  w_s(r)=\gamma_t(r),\text{ for }0\leq r< t\}. \end{align*}
$\dbM_{\gamma_t}$ is a subset of $\dbM$ containing all  paths
sharing the same path  $\gamma_t$ until time $t^-$.  Next, define $w_{t}^{-\e}\in \cM_{t-\e}$, the restriction of $w_t$ on $[0,t-\e)$, and $w_{t}^{+\e}\in \cM_{t+\e}$, the (continuous) extension of $w_t$ to $[0,t+\e)$, by
$w^{-\e}_{t}(s)= w_t(s),\q0\leq s< t-\e,$  and $ w^{+\e}_{t}(s)=\left\{\ba{ll}\ad w_t(s),\q0\leq s< t,\\
\ns\ad w_t(t^-),\q t\leq s<t+\e.\ea \right.$
We proceed with the following definitions.

\begin{definition} \label{defalphad}\rm A  measurable function  $h:\dbM\rightarrow\dbR$
is
\begin{itemize}
\item[{\rm (1)}] {\it path-right-continuous} (PR-continuous for short) if
$\lim_{\e\rightarrow 0^+} h(t+\e,w^{+\e}_t)=h(t,w_t)  \text{ for any $(t,w_t)\in \dbM$};$
\item[{\rm (2)}]	{\it path-left-continuous} (PL-continuous for short) if
$\lim_{\e\rightarrow 0^+} h(t-\e,w_t^{-\e})=h(t,w_t) \text{  for any $(t,w_t)\in \dbM$};$

\item[{\rm (3)}]{\it path-continuous} (P-continuous for short) if
it is PL- and PR-continuous.
Denote by $\sC^{0}(\dbM)$
the set of  P-continuous functions on $\dbM$;

\item[{\rm (4)}] {\it horizontally differentiable}, if there exists a $\partial^H h(t, w_t)$ such that
$$\partial^H h(t,w_t)=\lim_{\e\rightarrow 0^+}\frac{ h(t+\e,w_t^{+\e})-h(t,w_t)}{\e}=\lim_{\e\rightarrow 0^+}\frac{ h(t,w_t)-h(t-\e,w_t^{-\e})}{\e},$$
for any $(t,w_t)\in \dbM$;

\item[{\rm (5)}] {\it vertically differentiable}, if there exists $\{\partial_j^V h(t,w_t)\}_{j\in M}$ such that for any $(t,w_t)\in \dbM$,
$$\partial_j^V h(t,w_t)=\lim_{\e\rightarrow 0^+}h(t+\e,w_t\oplus j_{[t,t+\e)})=\lim_{\e\rightarrow 0^+}h(t,w^{-\e}_t\oplus j_{[t-\e,t)});$$

\item[{\rm (6)}] {\it $\alpha$-continuous}, if
$$\lim_{\e\rightarrow0^+}\dbE \big[h(t+\e,\alpha_{t+\e})\big|\a_t=w_t\big]=h(t,w_t)\text{  for any $(t,w_t)\in \dbM$};$$

\item[{\rm (7)}] {\it $\alpha$-differentiable},
if there exists a $\partial^\alpha h:\dbM\mapsto\dbR$ 
such that
\bel{ader}\partial^\alpha h(t,w_t)=\lim_{\e\rightarrow 0^+}\frac{\dbE [h(t+\e,\alpha_{t+\e})\big|\a_t=w_t]-h(t,w_t)}{\e}\text{  for any $(t,w_t)\in \dbM$};\eel

\item[{\rm (8)}]  A function $h:\dbM\mapsto\dbH$ is called  P-Lipschitz if there exists a uniform constant $\kappa_h>0$ such that  \bel{RsqrtLiP}
D_{\dbH}(h\big(t,w_r\oplus \gamma_{[r,t)}\big),h\big(t,\tilde w_r\oplus\gamma_{[r,t)}\big)\big)_{\dbH}\leq \kappa_h\sqrt{\sN(w_r)\sD(w_r,\tilde w_r)}\eel for  $\sD\big(w_r,\tilde w_r\big)<1$. Here $\dbH$ are some appropriate metric spaces with metric $D_{\dbH}(\cdot,\cdot).$ For example $\dbH$ can be $\dbR^n,$ $U$
 and $\sP_2(\dbR^n)$. Moreover, we always write $\kappa_h$ for the P-Lipschitz constant for the P-Lipschitz function $h$.
\end{itemize}
\end{definition}

The following lemma reveals
some general relationships among the above definitions.

\begin{lemma}\label{lemma-1} Let $h$ be a function defined on $\dbM$.
\begin{itemize}
\item[\rm (1)]  If $h$ is bounded, then
$h$ is $\alpha$-continuous if and only if $h$ is PR-continuous.

\item[\rm (2)] If $h$ is PR-continuous, then
$\partial^V_{w_t(t^-)}h(t,w_t)=h(t,w_t).$

\item [\rm(3)] If $h$ is  bounded and horizontally  and vertically differentiable  with
\bel{uniformrightcon}\lim_{\e\rightarrow 0^+}\sup_{\d\in(0,\e)}\Big |h(t+\e,w_t^{+\d}\oplus j_{[t+\d,t+\e)})-h(t+\e,w_t\oplus j_{[t,t+\e)})\Big|=0,\eel
and
\bel{uniformrightcon-1}\lim_{\e\rightarrow 0^+}\sup_{\d\in(0,\e)}\Big |h(t,w_t^{-\d}\oplus j_{[t-\delta,t)})-h(t,w^{-\e}_t\oplus j_{[t-\e,t)})\Big|=0,\eel
then $h$ is $\a$-differentiable with
$$\partial^\a h(t,w_t)=\partial^Hh(t,w_t)+\sum_{j\in M} q_{w_t(t^-),j}\partial_j^Vh(t,w_t).$$

\item[\rm(4)] For  $s>t$, let $\hat h(t,w_t):=\dbE [h(s,\a_{s})|\a_t=w_t].$
If $h$ is bounded,  P-continuous, and P-Lipschitz, then $\hat h(\cdot)$ is vertically and horizontally differentiable with P-continuous derivatives. As a consequence,  $\partial^\a \hat h(t,w_t)=0$ for $(t, w_t)\in\dbM$.

\item[\rm (5)]  If $h$ is horizontally differentiable with P-continuous derivative, then  there  exists $\d\in[0,\e)$ such that
\bel{mean-valuealpha}h(t+\e,w_t^{+\e})-h(t,w_t)=\e\partial^H h(t+\d,w_t^{+\d}).\eel

\end{itemize}	
\end{lemma}

\begin{proof}
(1)	Note that all the paths in $\dbM$ are  right-continuous. Denote by $N_\e$ the number jumps of $\alpha(\cdot)$ in $[t,t+\e)$. Given $\a_t=w_t\in \cM_t$, we have
$$
\ba{ll}\ad\!\!\!\dbE h(t+\e,\alpha_{t+\e})=\dbE\big[h(t+\e,\alpha_{t+\e})I(N_\e=0)\big]+ \dbE\big[ h(t+\e,\alpha_{t+\e})I(N_\e\geq 1)\big]\\
\ns\ad=h(t+\e,w_t^{+\e})\dbP(N_\e=0)+\Vert h\Vert_\infty O(\e)=h(t+\e,w_t^{+\e})+\Vert h\Vert_\infty O(\e).\ea
$$
Let $\e\rightarrow 0^+$,  the equivalence is clear.

(2) The claim is clear from
the definition of vertical derivatives.

(3) Set $\a_t=w_t$. On  $\{N_\e=1,\a((t+\e)^-)=j\neq w(t^-)\}$,
the possible values of $\a_{t+\e}$ are in the form of $w_t^{+\d}\oplus j_{[t+\d,t+\e)}\text{ for some }0\leq \d< \e.$
Note that
\bel{ito}\ba{ll}\ns\ad\dbE h(t+\e,\alpha_{t+\e})=\dbE\Big[h(t+\e,\alpha_{t+\e})\big[I(N_\e=0)+I(N_\e= 1)+I(N_\e\geq 2)\big]\Big]\\
\ns\ad=\!h(t\!+\!\e,w_t^{+\e})\dbP(N_\e=0)\!+\!\!\!\!\!\sum_{j\neq w_t(t^-)}^m\!\!\!\!\dbE\big[h(t+\e,\alpha_{t+\e})I(N_\e= 1,\a(t+\e^-)=j)\big]\!+\!O(\e^2)\\
\ns\ad=h(t+\e,w_t^{+\e})\big[1+q_{w_t(t^-),w_t(t^-)}\e+o(\e)\big]+ O(\e^2)\\
\ns\ad\q+\sum_{j\neq w_t(t^-)}\dbE\[[h(t+\e,\alpha_{t+\e})-\partial^V_jh(t,w_t)]I(N_\e= 1,\a(t+\e^-)=j)\]\\
\ns\ad\q+\sum_{j\neq w_t(t^-)}\partial^V_jh(t,w_t)\dbP(N_\e= 1,\a(t+\e^-)=j).\ea\eel
By \eqref{uniformrightcon}, it follows that
$\partial^\a h(t,w_t)
=\partial^H h(t,w_t) +\sum_{j\in M}q_{w_t(t^-),j}\partial^V_jh(t,w_t).$

(4)  Since $\a(\cdot)$ is a Markov process, we have
$$\ba{ll}\ns\ad \hat h (t+\e, w_t^{+\delta}\oplus j_{[t+\delta,t+\e)})=\dbE [h(s,w_t^{+\delta}\oplus j_{[t+\delta,t+\e)} \oplus \a_{[t+\e,s)})|\a(t+\e^-)=j]\\
\ns\ad  = \! \dbE [h(s,w_t \oplus j_{[t,t+\e)}\!\oplus \!\a_{[t+\e,s)})|\a(t+\e^-)=j] +o(1)\!=\!\hat h (t\!+\!\e, w_t\oplus j_{[t,t+\e)}) \!+\!o(1),\ea$$
where the $o(1)$ is uniformly depending on $\e$ but
independent of $\delta$. This verifies \eqref{uniformrightcon} and \eqref{uniformrightcon-1} can be verified similarly. Moreover, since  the probability of $\a$   jumps more than once in  $[t,t+\e)$ is $o(\e)$, we have
$$\ba{ll}\ds \lim_{\e\rightarrow 0^+}\hat h (t+\e, w_t\oplus j_{[t,t+\e)})=\lim_{\e\rightarrow 0^+}\dbE [h(s,w_t \oplus j_{[t,t+\e)}\oplus \a_{[t+\e,s)})|\a(t+\e^-)=j]\\
\ns\ds\q=\dbE [h(s,w_t \oplus \a_{[t,s)})|\a(t)=j].\ea$$
Similarly, we have
$$\ba{ll}\ds\lim_{\e\rightarrow 0^+}\hat h (t-\e, w^{-\e}_t\oplus j_{[t,t-\e)})=\lim_{\e\rightarrow 0^+}\!\dbE [h(s,w_{t-\e}  \oplus j_{[t-\e,t)}\oplus \a_{[t,s)})|\a(t^-)=j]\\
\ns\ds=\dbE [h(s,w_{t} \oplus \a_{[t,s)})|\a(t)\!=\!j].\ea$$
Therefore,  the vertical derivative of $\hat h(t,w_t)$ exists with $\partial^V \hat h(t,w_t)=  \dbE [h(s,w_t \oplus \a_{[t,s)})|\a(t)=j].$ Moreover, we can verify that   $\partial^V\hat h(t,w_t)$ is P-continuous.

Now, let us prove that $\hat h(t,w_t)$ is horizontally differentiable.
Let $\a_t=w_t$ and $s>t$. Note that if $\a$ has only one jump in $[t,t+\e)$ and $\a(t+\e^-)=j$, we have
$
\sD(w_t\oplus\a_{[t,t+\e)}, w_t\oplus j_{[t,t+\e)})\leq \e.
$
Then, it follows that
$$\ba{ll}\ns\ad\hat h(t+\e,w_t^{+\e})-\hat h(t,w_t)=\hat h(t+\e,w_t^{+\e})-\dbE [h(s,\a_{s})]\\
\ns\ad=\hat h(t+\e,w_t^{+\e})\!-\!\dbE[ h(s,\a_{s})I(N_\e=0)]\!-\!\dbE [h(s,\a_{s})I(N_\e=1)\!+\! h(s,\a_{s})I(N_\e=2)]\\
\ns\ad=\hat h(t+\e,w_t^{+\e})-\hat h(t+\e,w_t^{+\e})(1+q_{ii}\e+O(\e^2))+\Vert h\Vert_\infty O(\e^2)\\
\ns\ad\q-\sum_{j\neq i}\hat  h(t+\e,w_t\oplus j_{[t,t+\e)})\dbE I(N_\e=1,\a(t+\e^-)=j)\\
\ns\ad\q+\sum_{j\neq i}\dbE[   h(s,w_t\oplus j|_{[t,t+\e)}\oplus \a_{[t+\e,s)})-h(s,\a_s)]I(N_\e=1,\a(t+\e^-)=j)\\
\ns\ad=-\e\sum_{j\neq w_t(t^-)} q_{w_t(t^-),j}\partial^V_j\hat h(t, w_t)+o(\e).\ea
$$
Similarly,
$\hat h(t,w_t)-h(t-\e,w_t^{-\e})=-\e\sum_{j\in M} q_{w^{-\e}_t(t-\e^-),j}\partial^V_j\hat h(t-\e, w_t^{-\e})+o(\e).$
By the P-continuity of $\partial^V_j\hat h(t, w_{t})$, we have $\hat h(t,w_t)$ is horizontally differentiable with
$\partial ^H\hat h(t,w_t)=-\sum_{j\in M} q_{w_t(t^-)j}\partial^V_j\hat h(t, w_t).$
Moreover,
$\partial ^H\hat h(t,w_t)$ is also P-continuous. Consequently,  $\partial^\a \hat h(t,w_t)=0$.

(5) As long as $ h(t,w_t)$ is horizontally differentiable with P-continuous derivative, \eqref{mean-valuealpha} holds
by the mean-value theorem.
\end{proof}

We emphasize that the P-Lipschitz property is new in the literature, which is used to guarantee the $\a$-differentiability of the solutions of  path-dependent HJBs. In fact, if the optimal control problem does not involve mean-field interactions, then no path-dependence is needed,
and
its value function takes the form of $V(t,w_t,x)=V(t,w_t(t^-),x)$,
which is naturally P-Lipschitz with P-Lipschitz constant 0. Then the study on P-Lipschitz property  is not needed for such special scenario. Because we need to consider the path-dependence in our paper,   we need to generalize the calculus on the path-space from   the classical  theories for switching diffusions without mean-field terms.
With the path-space of $\a$, we will focus on the HJB-FP approach for our problems in the next two sections.

\section{FP-Step}\label{sec:fp}
This section is devoted to the FP-step. 
 Let $L^2_{\dbF}(\Omega; C([0,T],\dbR^n))$ be the  Banach space of $\dbF$-progressively measurable processes with continuous paths equipped with the sup-norm
$\Vert X\Vert^2_{\sL^2}:=\dbE\sup_{0\leq t\leq T}|X(t)|^2.$ First, we present a proposition concerning the existence and uniqueness of a solution to SDE \eqref{SDE-1} for $u(\cdot)\in\cU[0,T]$. The proof is based on a standard Picard's iteration and thus is omitted here.

 \begin{proposition}\label{exsolutionsde}   Suppose Assumption \ref{A-3} holds. 
For any $u(\cdot)\in \cU[0,T]$,  there exists a unique solution $X\in L^2_{\dbF}(\Omega, C([0,T],\dbR^n))$ to \eqref{SDE-1}.


\end{proposition}





To construct the mapping $\zeta$, we need to  restrict $u$ in a feedback strategy instead. 
Let  $\sU$ be the set of {\it Lipschitz feedback strategies} defined by
\bel{ads}\ba{ll}\ns\ds \sU:=\bigcup_{\kappa\geq 0}\sU_\kappa\text{ and }\\
\ds\sU_\kappa\!\!:=\!\!\Big\{u\!:\!\dbM\!\times\! \dbR^n \! \mapsto \!U\!\Big| u(\cdot) \text{ is P-continuous, uniformly P-Lipschitz with constant $\kappa$,}\\
\ns\ds\qq\qq\text{ and } |u(t,w_t,0)|\leq \kappa,~|u(t, w_t,x)-u(t, w_t,y)|\leq \kappa |x-y|
\Big\}.\ea\eel
 Given $u\in\sU$, inspired by \cite{NguyenYinHoang 2020}, let us construct $\zeta$  as follows.

For any $(t,w_t)\in\dbM$,   consider the following SDE on $[0,t]$,
\bel{givenpathZ} dY(s)=\Bb(s, w(s^-),Y(s),\text{law}(Y(s)); u(s,w_s,Y(s)))ds+
\BBsi(s,Y(s))dB(s),\eel
with $\text{law}(Y(0))=\mu_0$. Because $(t,w_t)\in\dbM$ is given, one can solve
the above SDE until time $t$.
Set
$\zeta_{\mu_0}^u(t,w_t):=\text{law} (Y(t))$
to emphasize the dependence of $\zeta$ on $u$ and $\mu_0$.
The following proposition is a direct consequence of Theorem 4.6  in \cite{NguyenYinHoang 2020}, which guarantees  that $\zeta_{\mu_0}^u(t,\a_t)$ is the correct feedback function for the conditional distribution $\eta(t)$.

\begin{proposition}\label{prop-1} Suppose Assumption \ref{A-3} holds.
Let  $X^u$ be the solution to \eqref{SDE-1} under strategy $u\in\sU$ and recall that
$\eta^u=\{\eta^u(t)=\law(X^u(t)|\cF_{t}^\a):0\leq t\leq T\}.$
Then,
$$\dbP\big(\eta^u(\cdot)=\zeta^u_{\mu_0}(\cdot,\alpha_\cdot)\text{ in } D([0,T],\cP_2(\dbR^n)) \big)=1.$$
Here $D([0,T],\cP_2(\dbR^n))$ denotes the space of $\cP_2(\dbR^n)$-valued c\'ad\'ag curves on $[0,T]$   endowed with the Skorohod metric.
\end{proposition}

With the exact form of the feedback function $\zeta$ given a $u\in\sU$,  let us present the following lemma on P-Lipschitz property of $\zeta$ which will be used in the future.
\begin{lemma}\label{uzeta}  Suppose Assumption \ref{A-3} holds and  let $u\in\sU_{L_0}$.
Then $\zeta_{\mu_0}^u(t,w_t)$ is P-Lipschitz with a constant $\sqrt{\beta_0(T)}(1+ L_0)$.
Moreover,  for any $u,\tilde u\in\sU$,
\bel{zetauniform}
\sup_{(t,w_t)\in\dbM} \Bw\big(\zeta_{\mu_0}^u(t,w_t),\zeta_{\mu_0}^{\tilde u}(t, w_t)\big)\leq \sqrt{\beta_1(T)}\sup_{(t,w_t,x)\in\dbM \times\dbR^n}\big|u(t,w_t,x)-\tilde u(t,w_t,x)\big|.
\eel
Here the two constants $\beta_0(T)$ and $\beta_1(T) $ are independent of $L_0$ and are  small when $T$ is small.
\end{lemma}

\begin{proof}   Let $Z(\cdot)$ and $\widetilde Z(\cdot)$ be the solutions to \eqref{givenpathZ} under the same strategy $u$ and the same initial
distribution $\mu_0$, but with different paths $w_t= w_r\oplus\gamma_{[r,t)}$ and $\tilde w_t=\tilde w_r\oplus\gamma_{[r,t)}$, respectively.
Then, for $s\in[0,t]$, It\^o's formula implies that
\begin{align}\label{ZZZ}& \frac d {ds} \dbE\big|Z(s)-\tilde Z(s)\big|^2=\big|\BBsi(s,Z(s))-\BBsi(s,\tilde Z(s))\big|^2\\
& \q+2\dbE\big\lan\! Z(s)\!-\!\tilde Z(s),\Bb\big(\!s,Z(s),w(s^-),\law\!\big(Z(s)\big),u\big(s,w_s,Z(s)\big)\!\big)\nonumber\\
&\qq\qq\qq\qq-\Bb\big(\!s,\tilde Z(s),\tilde w(s^-),\law\!\big(\tilde Z(s)\big),u\big(s,\tilde w_s,\tilde Z(s)\big)\!\big)\big\ran\nonumber\\
&\leq K(1+L_0)\big[\dbE|Z(s)-\tilde Z(s)|^2+\sN(w_r)\sD\big(w_r,\tilde w_r\big)+I\big(w(s)\neq \tilde w(s),0\leq s\leq r\big)\big].\nonumber
\end{align}
  By \eqref{SNSD},  Gronwall's inequality yields $\dbE|Z(t)-\wt Z(t)|^2\leq Kt(1+L^2_0)\sN(w_r)\sD\big(w_r,\tilde w_r\big)$. Note that  $Z(t)$ and $\wt Z(t)$ have distribution $\zeta_{\mu_0}^u(t,w_t)$ and $\zeta_{\mu_0}^u(t,\tilde w_t)$ respectively, by the definition of Wasserstein-2 metric, we know that $\zeta_{\mu_0}^u(t,w_t)$ is P-Lipschitz with constant $KTL_0$.

To prove \eqref{zetauniform}, let $Z(t)$ and $\widetilde Z(t)$ be the solutions to \eqref{givenpathZ} with the same initial  distribution
 and
path $w$ but with different strategies $u$ and $\tilde u$, resp.
 Similar to \eqref{ZZZ}, one has $$\frac d{ds}\dbE\big|Z(s)-\tilde Z(s)\big|^2\leq K\dbE\big|Z(s)-\tilde Z(s)\big|^2+ K\sup_{(t,w_t,x)\in\dbM \times\dbR^n}\big|u(t,w_t,x)-\tilde u(t,w_t,x)\big|^2.
$$
 Using  Gronwall's inequality again, we have
\eqref{zetauniform} similarly.
\end{proof}



At the end of this section,  we denote by $\Upsilon$ the set of all such $\zeta_{\mu_0}^u$ for any $u\in\sU$ and $\mu_0\in\sP_2(\dbR^n)$.  To show the dependence on the initial $\mu_0$, we write $\Upsilon_{\mu_0}:=\big\{\zeta\in\Upsilon:\zeta(0,w_0)=\mu_0\big\}$. We also equip  $\Upsilon$ with the supremum metric,
$
d_\Upsilon(\zeta,\tilde\zeta):=\sup_{(s,w_s)\in\dbM} d\big(\zeta(s,w_s),\tilde\zeta(s,w_s)\big).
$
Now, we can complete the  FP-step by defining a map $\cT_1^{\mu_0}:\sU\mapsto \Upsilon_{\mu_0}$ by \bel{T1}[\cT^{\mu_0}_1u](s,w_s):=\zeta_{\mu_0}^u(s,w_s)\text{ for all }(s,w_s)\in\dbM.\eel
We will work with the HJB-step in the next section.

\section{HJB-Step} \label{sec: HJB-step}
In this section,  we
focus on the HJB-step. Replacing $\eta(t)$ by the feedback form $\zeta(t,\a_t)$,  the HJB-step is equivalent to solving a   path-dependent  but distribution-independent control problem under general discounting cost. As mentioned, we will tackle the control problem by studying the equilibria of HJBs.  
We note that compared to the path-dependence theory using
viscosity solutions developed in the
literature such as \cite{EKren2016a,EKren2016b,Cosso2022}, the solutions constructed here are in the classical sense due to the special structure   of the path space. Such strong solutions are necessary for us to work with the path-dependent equilibrium HJBs \eqref{PHJE-10000}.
Write \begin{align*}&b(t,w_t,y;v)=\Bb(t,y,w_t(t^-),\zeta(t,w_t);v),~ \sigma(t,x)=\BBsi(t,x),\\&g(\t;t,w_t,y;v)=\Bg(\t;t,y,w_t(t^-),\zeta(t,w_t);v),~
h(\t;w,y)=\Bh(\t;w(T^-),y,\zeta(t,w)).\end{align*}
Here $b$, $g$, and $h$ are dependent on $\zeta$ while we omit this dependence in the notations for simplicity.

Given $\zeta$, the SDE in \eqref{SDE-1} reduces to
\bel{PDSDE}dX(t)=b(t,\alpha_{t},X(t);u(t))dt+\sigma(t,X(t)) dB(t).\eel
 The general  discounting cost functional  becomes
\begin{align}\label{PDLLLLL}&V(t, w_t,x;u(\cdot))=J(t;t, w_t,x;u(\cdot)),\\&J(\tau;t, w_t,x;u(\cdot))=\dbE_{t,w_t,x}\(\int_t^T g(\tau;s,\alpha_{s},X(s);u(s))dt+h(\tau;\alpha_{T},X(T))\).\nonumber\end{align}
Our goal is to find a closed-loop  strategy for the control problem subject to \eqref{PDSDE} under the cost \eqref{PDLLLLL}.
Due to the general discounting structure,  the control problem is time-inconsistent where no optimal control exists.
Therefore, we turn to find a local optimal {\it equilibrium strategy} inspired by \cite{Yong2012} where the key step lies in solving an equilibrium for a systems of HJBs. In our case,   the equilibrium HJBs turn out to be path-dependent, leading to
\begin{align}\label{PDEHJB-00}\begin{cases}&\!\!\!\!(\partial^\a\!+\!\cA)\Th(\tau;t,w_t,x) +\wt \cH(\t;t,w_t,x, D_x\Theta(\t;t,w_t,x),D_x\Theta(t;t,w_t,x))=0,\\
&\!\!\!\!\Theta(\tau;T,w,x)=h(\tau;w,x),\end{cases}\end{align}
where
$\wt \cH(\t;t,w_t,x, p,q)=\cH(\tau;t, w_t(t^-), x,\zeta(t,w_t),p,q).$
Here, $D_x$ and $D_x^2$ are the gradient and Hessian matrix with respect
to $x$ and recall that  $\partial^\a$  has been defined in \eqref{ader}. Compared to the classical HJB, the additional time factor $\tau$ is introduced due to the general discounting structure.  Our goal is to show (i) there exists a unique classical solution to \eqref{PDEHJB-00}, which is regular enough such that a feedback strategy can be identified, (ii)  the feedback strategy is local optimal in some appropriate sense.
The two assertions are to be
studied in the following subsections separately. Before moving on, we present the following proposition  whose proof is straightforward by Assumption \ref{A-3} and thus is omitted.

\begin{proposition}\label{defvarphi0} Suppose Assumption \ref{A-3} holds and $\zeta$ is P-continuous and  P-Lipschitz.

{\rm (1)} 	 The mappings $(t,w_t,x,v)\mapsto b(t,w_t,x,v)$ and  $(\t;t,w_t,x,v)$ $\mapsto g(\t;t,w_t,x,v) $ are P-continuous and continuous with respect to $(\t,t,x,v)$;   $h(\t;w,x) $ is continuous in $\t$. It satisfies that
$$\left\{\ba{ll} \ad|b(t,w_t,x,v)|+|g(\t;t,w_t,0,0)|+
|h(\t;w,0)|\leq K,\\
\ns\ad|b(t,w_t,x_1,v_1)-b(t,w_t,x_2,v_2)|\leq K(|x_1-x_2|+|v_1-v_2|),\\
\ns \ad|g(\t;t,w_t,x_1,v_1)-g(\t;t,w_t,x_2,v_2)\leq K(|x_1-x_2|+|v_1-v_2|),\\
\ns \ad|D_xh(\t;w,x)|+|D_x^2h(\t;w,x)|\leq K,\\
\ns\ad b, D_xh\text{ and } D_x g\text{ are P-Lipschitz with $\kappa_b,\kappa_{D_xh},\kappa_{D_xg}\leq K \kappa_\zeta$}
\ea\right.$$
 {\rm (2)}
 The Hamiltonian  $\wt\cH(\t;t,w_t,x,p,q)$  is P-continuous with $|\wt\cH(\cdot)|\leq K(1+|p|+|x|),~|D_q\wt\cH(\cdot)|\leq K(1+|p|),~ |D_p\wt\cH(\cdot)|\leq K$ and
$$|\wt\cH(\t;t,w_r\oplus\gamma_{[r,t)},x,p,q)-\wt\cH(\t;t,\tilde w_r\oplus\gamma_{[r,t)},x,p,q)|\! \leq \!K\kappa_{\zeta}(1+|p|)\sqrt{\sN(w_r)\sD(w_r,
\tilde w_r)},$$
if $\sD(w_r,\tilde w_r)<1$.

\end{proposition}

\subsection{Path-Dependent  PDEs}
To work on path-dependent equilibrium HJBs
\eqref{PDEHJB-00},  we will  adopt the classical fundamental solution method. To better illustrate the details, we start with a path-dependent PDE first instead of the complex equilibrium HJBs
\eqref{PDEHJB-00}.

 For some appropriate functions $b$, $g$, and $h$, we consider the following path-dependent PDE,
\bel{PDPDE222} \left\{\ba{ll}\ns \ad (\partial^\alpha+\cA)\Th(t,w_t,x)+\lan b(t,w_t,x), D_x\Th(t,w_t,x)\ran +g(t,w_t,x)=0, \\
\ns \ad\Th(T,w,x)=h(w,x).
\ea\right.\eel
In this subsection only, $b$, $g$, and $h$ are not necessarily the same as those in the previous sections. The key idea of this subsection lies in applying the fundamental solution method for parabolic PDEs. Then we will extend the results to the equilibrium HJBs in our paper.

Due to $\partial^\a$ and $\cA$,  the fundamental solution consists of two parts here. Define
\bel{heatkernel}\Psi_1(t,x;s,y):=\frac1{(4\pi (s-t))^{\frac n2}(\text{det}[a(s,y)])^{\frac12}}\exp\Big\{-\frac{\lan x-y, a^{-1}(s,y)(x-y)\ran}{4(s-t)}\Big\},\eel
which is the fundamental solution to
$\partial _s \pi(s,x)+\cA \pi(s,x)=0.$
Straightforward calculation yields that
$$\left\{\ba{ll}\ds|\Psi_1(t,x;s,y)|\leq K(s-t)^{-\frac n2}\exp\Big\{-\l_0\frac{|x-y|^2}{4(s-t)}\Big\},\\
\ns\ds |D_y\Psi_1(t,x;s,y)|\leq K(s-t)^{-\frac {n+1}2}\exp\Big\{-\l_0\frac{|x-y|^2}{4(s-t)}\Big\},\ea\right.$$
and
\bel{Dypsi}D_y\Psi_1(t,x;s,y)=-D_x\Psi_1(t,x;s,y)+\Psi_1(t,x;s,y)\rho(t,x;s,y),\eel
where
$$\left\{\ba{ll}\ds
\rho(t,x;s,y)=-\frac{D_y(\text{det}[a(s,y)])}{2\text{det}[a(s,y)]}-\frac{\lan [ a^{-1}(s,y)]_y(x-y),x-y\ran}{4(s-t)},\\ [2mm]
\ns\ds\lan  [a^{-1}(s,y)]_y(x-y),x-y\ran=\begin{pmatrix}[a^{-1}(s,y)]_{y_1}(x-y),x-y\\\vdots\\ [a^{-1}(s,y)]_{y_n}(x-y),x-y\end{pmatrix}. \ea\right.$$
It can be checked  that under our assumption, one has
\bel{rhoesti}|\rho(t,x;s,y)|+|D_y\rho(t,x;s,y)|\leq K\(1+\frac{|x-y|^2}{s-t}\).\eel

Let us introduce the second part of the fundamental solution from  the path process. Since $\a(\cdot)$ is a homogeneous Markov chain, we can define
$\Phi_2(i,d\gamma_t):=\dbP(\a_t\in d\gamma_t|\a(0)=i).$
It can be seen that $\Phi_2(i,d\gamma_t)$ is supported in $\gamma_t(0)=i$ only.
Then  the transition law of $\a_\cdot$ can be written as
$\Psi_2(t,w_t;s,d\tilde w_s):=\dbP(\a_{s}\in d\tilde w_s|\a_{t}=w_t)
=\Phi_2(w_t(t^-),d\gamma_{[t,s)})$
where $\gamma_{[t,s)}(r)=\tilde w_s(r)$ for $t\leq r< s$.

Write
$$\Psi(t,w_t,x;s,d\tilde w_s,y):=\Psi_1(t,x;s,y)\Psi_2(t,w_t;s,d\tilde w_s).$$
The	function		
$\Psi$  plays a role as a fundamental solution in our path-dependent PDE. Note that $\Psi_2(t,w_t;s,\cdot)$ has support in $\{\tilde w_s\in\cM_s:\tilde w_s(r)=w_t(r), r\in[0,t)\}$ only. Before  presenting our results, we first introduce
some
notation of the spaces for the solution to our path-dependent PDE.

Define $\sC^{0,0}(\dbM\times\dbR^n:\dbR^m)$ by the set of all $\dbR^m$-valued functions $\th(t,w_t,x)$ which are P-continuous for  each fixed $x$ and are continuous with respect to $x$ for each $w_t\in \cM_t$.
Define the sup-norm on $\sC^{0,0}(\dbM\times\dbR^n:\dbR^m)$ by
$\Vert\theta\Vert_{\sC^{0,0}}:=\sup_{\dbM\times\dbR^n}|\theta|.$
Note that
$\dbM_{w_t}$ is a subset of $\dbM$, $\sC^{0,0}(\dbM_{w_t}\times\dbR^n:\dbR^m)$ can be  defined similarly. For $\theta\in \sC^{0,0}(\dbM_{w_t}\times\dbR^n:\dbR^n)$, the sup-norm  is defined by $\Vert\theta\Vert_{\sC_{w_t}^{0,0}}
:=\sup_{\dbM_{w_t}\times\dbR^n}|\theta|.
$
Note that $\sC^{0,0}(\dbM\times\dbR^n:\dbR^m)$ and  $\sC^{0,0}(\dbM_{w_t}\times\dbR^n:\dbR^m)$ are  complete  under the respective norms.
Now, we can solve \eqref{PDPDE222}.

\begin{proposition}\label{fundamental}
Let $b(t,w_t,x)$ and    $g(t,w_t,x)$ be P-continuous for each  $x$ and are continuously differentiable with respect $x$. Suppose $h(w,x)$ is twice continuously differentiable in $x$. Further suppose that $ g_x(t,w_t,x)$ is P-continuous and
\bel{conditionfg}\left\{\ba{ll}
\ad|b(t,w_t,x)|+|b_x(t,w_t,x)|+|g(t,w_t,0)|+|h(w,0)|\leq K,\\
\ns\ad |g_x(t,w_t,x)|+|h_x(w,x)|+|h_{xx}(w,x)|\leq K,\\
 \ns\ad b(t,w_t,x), g(t,w_t,x), h_x(w,x) \text{ are P-Lipschitz uniformly in $x$}. 
\ea\right.\eel
Then the path-dependent PDE \eqref{PDPDE222}
admits a unique solution $\Theta$ in the classical sense in the class of functions satisfying
$|\Th(t,w_t,x)|\leq Ce^{\delta x^2}\text{ for any }\delta>0.$  The solution
has  the following representation
\bel{PDEX}\ba{ll}\ns\ad\Th(t,w_t,x)=\int_{\cM\times \dbR^n}h(\tilde w, y)\Psi(t,w_t,x;T,d\tilde w,y)dy\\
\ns\ad+\!\int_t^T\!\int_{\cM_s\times \dbR^n}\!\!\(\lan b(s,\tilde w_s,y), D_y\Th(s,\tilde w_s,y)\ran+g(s,\tilde w_s,y)\)\Psi(t,w_t,x;s,d\tilde w_s,y)dyds.\ea\eel
 Moreover, $D_x\Th(t,w_t,x),~D_x^2\Th(t,w_t,x),~\partial^\a \Th(t,w_t,x)$ are P-continuous for each $x$ and is continuous  w.r.t. $x$ for each $(t,w_t)\in\dbM_t$. In addition, $D_x\Th(t,w_t,x)$ is P-Lipschtiz  with  Lipschitz constant $\beta_2(T)(1+\kappa_\zeta)$ for some constant $\beta_2(T)>0$.
 Here the  solution in the classical sense means that $\Th(t,w_t,x)$ is
 twice
 differentiable in $x$, $\a$-differentiable in $(t,w_t)$ (see Definition \ref{defalphad}) and satisfies \eqref{PDPDE222}.
\end{proposition}

\begin{proof}
The proof is based on a fixed point theory.
We first prove there exists a solution on $[T-\delta,T]$ for some positive constant $\delta$.
Then similar idea can be extended to the whole time interval $[0,T]$. Our proof is divided into several steps.


(1) {\it Construction of the map $\Pi$ on $[T-\d,T]$}. We fix a $(T-\d,\gamma_{T-\d})\in\dbM$ and take  a $\th\in\sC^{0,0}(\dbM_{\gamma_{T-\d}}\times\dbR^n:\dbR^n)$. We will define a map $\Pi:\sC^{0,0}(\dbM_{\gamma_{T-\d}}\times\dbR^n:\dbR^n)\mapsto\sC^{0,0}(\dbM_{\gamma_{T-\d}}\times\dbR^n:\dbR^n)$ such that $\Pi$ admits a fixed point.

For any $(t,w_t,x)\in \dbM_{\gamma_{T-\d}}\times\dbR^n$, define
\begin{equation}\label{thtoTh}\Th(t,w_t,x):=\int_{\dbR^n}\Psi_1(t,x;T,y)\hat h(t,w_t,y)dy+\int_t^T\int_{\dbR^n}\Psi_1(t,x;s,y)\hat g(t,w_t,s,y)dyds.\end{equation}
where \begin{align*}
&\hat h(t,w_t,y):=\int_{\cM} h(\tilde w, y)\Psi_2(t,w_t;T,d\tilde w)=\int_{\cM[t,T)} \!\!\!\!h(w_t\oplus \gamma_{[t,T)},y)\Phi_2(w_t(t^-),d\gamma_{[t,T)});
\\
&\hat g(t,w_t,s,y):=\int_{\cM_s}\(\lan b(s, \tilde w_s,y), \theta(s, \tilde w_s,y)\ran+g(s, \tilde w_s,y)\)\Psi_2(t,w_t;s, d\tilde w_s)\\
&=\int_{\cM[t,s)}\!\!\!\(\lan b(s,\! w_t\!\oplus\! \gamma_{[t,s)},\!y)\!,\! \theta(s, w_t\oplus \gamma_{[t,s)},y)\ran\!+\!g(s, w_t\oplus \gamma_{[t,s)},y)\)\Phi_2(w_t(t^-)\!,d\gamma_{[t,s)}).
\end{align*}
By the P-continuity of $b,\th,g$, and $|D_y h|<K$, it can be seen that $\hat h(t,w_t,y)$ is uniformly Lipschitz in $y$ and
$\hat g(t,w_t,s,y)$ is continuous in $(s,y)$.  By \cite[Theorem 12, p25]{Friedman}, $\Th(t,w_t,x)$ is second order differentiable w.r.t $x$. Using integration by parts, we have \begin{align}\label{repDTh}&\ds D_x\Th(t,w_t,x)=\int_{\dbR^n}D_x\Psi_1(t,x;T,y)\[\int_{\cM}h(\tilde w,y)\Psi_2(t,w_t;T,d\tilde w) \]dy\\\nonumber
&\ds +\!\int_t^T\!\!\!\int_{\dbR^n}D_x\Psi_1(t,x;s,y) \int_{\cM_s}\!\!\![\lan b(s,\tilde w_s,y), \theta(s,\tilde w_s,y)\ran\!+\!g(s,\tilde w_s,y)]\Psi_2(t,w_t;\!s,d\tilde w_s) dy ds\\\nonumber
&\ds=\int_{\dbR^n}\!\!\![\Psi_1(t,x;T,y)\rho(t,x;T,y)-D_y\Psi_1(t,x;T,y)]\[\int_{\cM}h(\tilde w,y)\Psi_2(t,w_t;\!T,d\tilde w) \]dy\\\nonumber
&\ds+\!\int_t^T\!\!\!\int_{\dbR^n}D_x\Psi_1(t,x;s,y) \int_{\cM_s}\!\!\![\lan b(s,\tilde w_s,y), \theta(s,\tilde w_s,y)\ran\!+\!g(s,\tilde w_s,y)]\Psi_2(t,w_t;\!s,d\tilde w_s) dy ds\\\nonumber
&\ds=\int_{\dbR^n}\!\!\!\Psi_1(t,x;T,y)\[\int_{\cM}[D_y +\rho(t,x;T,y)] h(\wt w,y)\Psi_2(t,w_t;T,d\tilde w)\] dy\\\nonumber
&\ds +\!\int_t^T\!\!\!\int_{\dbR^n}D_x\Psi_1(t,x;s,y) \int_{\cM_s}\!\!\![\lan b(s,\tilde w_s,y), \theta(s,\tilde w_s,y)\ran\!+\!g(s,\tilde w_s,y)]\Psi_2(t,w_t;\!s,d\tilde w_s)dyds.
\end{align}	
By \eqref{Dypsi} and  the P-continuity of $b,\th,g$, we see that $D_x\Th(t,w_t,x)$ is P-continuous  with
\begin{equation}\label{lipthTh}\sup_{(t,w_t,x)\in \dbM_{\gamma_{T-\d}}\times\dbR^n}|D_x\Th(t,w_t,x)|\leq \ell_0(1+\sqrt{\d})\sup_{(t,w_t,x)\in \dbM_{\gamma_{T-\d}}\times\dbR^n}|\th(t,w_t,x)|<\infty.\end{equation}
Here $\ell_0$ is a
constant independent of $\th$ and $t$.
Now we can define the map $\Pi$ by
$\Pi[\th]=D_x\Th$.

(2) {\it Fixed-point of $\Pi$}. By \eqref{repDTh}, we have a constant $\ell_1$ depending on $b$ and $\sigma$ only such that $$\ba{ll}\ns\ad |\Pi[\th_1](t,w_t,x)-\Pi[\th_2](t, w_t,x)|\\
\ns\ad\leq \int_t^T\int_{\dbR^n}D_x\Psi_1(t,x;s,y) \[\int_{\cM_s}|\lan b, \theta_1- \theta_2\ran|(s,\tilde w_s,y)\Psi_2(t,w_t;s,\tilde w_s) \]dy ds\\
\ns\ad\leq \ell_0\sqrt \d\Vert\theta_1-\theta_2\Vert_{\sC_{\gamma_{T-\delta}}^{0,0}}.\ea$$
Take $\d=\min(\ell_0^{-2}/2,\ell_1^{-2}/2)$, we see that $\Pi$ is a contraction on $\sC^{0,0}(\dbM_{\gamma_{T-\d}}\times\dbR^n:\dbR^n)$  which admits a unique fixed point $\th^*$. Write $\Th^*$ be the function $\Th$ in \eqref{thtoTh} by taking $\th=\th^*$ and $\Th^*$ is the solution to \eqref{PDEX} for $(t,w_t,x)\in \dbM_{\gamma_{T-\d}}\times\dbR^n.$   Note that $\d$ is independent of the choice  $\gamma_{T-\d}$ and we can take arbitrary $\gamma_{T-\d}\in \cM_{T-\d}.$ Then we construct a solution $\Th^*$ to \eqref{PDEX} for $t\geq T-\d$. By
\eqref{lipthTh}, we have $|D_x\Th(T-\d,w_{T-\d},x)|\leq \ell_0(1-\ell_0\sqrt\d).$

(3) {\it Extension to the full time horizon $[0,T]$}. Repeat the above process on $[T-2\d,T-\d]$ with $\d=\ell_0^2/2$ and use $\Th(T-\d, w_{T-\d},x)$ as a terminal condition. Note that $\d=\min(\ell_0^{-2}/2,\ell_1^{-2}/2)$ is an absolute constant independent of time, we can get a solution on $[T-2\d,T-\d]$ similarly. Hence, we are allowed to repeat the above process to $[T-2\delta,T-\delta]$, $[T-3\delta,T-2\delta]$, $\ldots$, until reaching the initial time 0. Then \eqref{PDPDE222}  admits a unique  solution  $\Theta^*$ on $[0, T]$ with $\th^*=D\Th^*\in \sC^{0,0}(\dbM\times\dbR^n:\dbR^n)$.

(4)  {\it Verification of classical solution}. Now we want to verify $\Th^*$ is the classical solution to path-dependent PDE  \eqref{PDPDE222}. First, we know that $\Th^*$ is 
twice
differentiable w.r.t. $x$ with
\begin{align}\label{lapaceTh}&\cA\Th^*(t,w_t,x)=\frac12\text{Trace}[a(t,x) D_x^2\Th^*(t,w_t,x)]\\
&=\frac12\int_{\dbR^n}\text{Trace}[a(t,x)D^2_x\Psi_1(t,x;T,y)]\hat h(t,w_t,y)dy\nonumber\\
&\q+\frac12\int_t^T\int_{\dbR^n}\text{Trace}[a(t,x)D^2_x\Psi_1(t,x;T,y)] \hat g(t,w_t,s,y)dy ds\nonumber\\
&=-\int_{\dbR^n}\partial_t\Psi_1(t,x;s,y)\hat h(t,w_t,y)dy-\int_t^T\int_{\dbR^n}\partial_t\Psi_1(t,x;s,y) \hat g(t,w_t,s,y)dy ds.\nonumber
\end{align}
Then we want to show that $\Th^*$ is $\a$-differentiable. By  Lemma \ref{lemma-1} and \eqref{thtoTh}, it suffices to show $\hat h(t,w_t,x)$ and $\hat g(t,w_t,s,y)$ are P-Lipschitz in $(t,w_t).$ The P-Lipschitz property of $\hat h$ is true by that of $h$ directly. By the definition of  $\hat g(t,w_t,s,y)$, we need only  show $\th^*$ is  P-Lipschitz.

For $\sD(w_r,\tilde w_r)<1$, we know that $w_r(r^-)=\tilde w_r(r^-)$. For $s\geq t$, write
$$
\kappa_{\th^*}(s)=\sup_{\gamma_{[r,s')}\in\cM[r,s'),s'\geq s}\sup_{x\in\dbR^n}\frac{|\th^*(t,w_r\oplus  \gamma_{[r,s')},x) -\th^*(t,\tilde w_r\oplus  \gamma_{[r,s')},x) |}{\sqrt{\sN(w_r)\sD(w_r,\tilde w_r)}}.
$$
Note that
for $t\geq r$,  \begin{align}\label{lipTh2}&|D_x\Th^*(t,w_r\oplus \gamma_{[r,t)},x)-D_x\Th^*(t,\tilde w_r\oplus \gamma_{[r,t)},x)|\\
&\leq \int_{\cM\times \dbR^n}|(D_y+\rho(t,x;T,y)) h(w_r\oplus \gamma_{[r,t)}\oplus\g_{[t,T)}, y)\nonumber\\
&\qq-\!(D_y\!+\!\rho(t,x;T,y))h(\tilde w_r\oplus \gamma_{[r,t)}\oplus\g_{[t,T)}, y)|\Psi_1(t,x;T,y)dy\Phi_2(\gamma_{[r,t)}(t^-),d\g_{[t,T)})\nonumber\\
&\q+\int_t^T\int_{\cM_s\times \dbR^n}|\lan b,\th^*\ran (s,w_r\oplus \gamma_{[r,t)}\oplus\g_{[t,s)}, y)-\lan b,\th^*\ran(s,\tilde w_r\oplus \gamma_{[r,t)}\oplus\g_{[t,s)}, y)|\nonumber\\
&\qq\qq\times|D_x\Psi_1(t,x;s,y)|dy\Phi_2(\gamma_{[r,t)}(t^-),d\g_{[t,s)})ds\nonumber\\
&\q+\int_t^T\int_{\cM_s\times \dbR^n}|g(s,w_r\oplus \gamma_{[r,t)}\oplus\g_{[t,s)}, y)-g(s,\tilde w_r\oplus \gamma_{[r,t)}\oplus\g_{[t,s)}, y)|\nonumber\\
&\qq\qq\times|D_x\Psi_1(t,x;s,y)|dy\Phi_2(\gamma_{[r,t)}(t^-),d\g_{[t,s)})ds\nonumber\\
&\leq \ell_2(1+\kappa_\zeta)\(1+\int_t^T\frac{1}{\sqrt{T-s}}\kappa_{\th^*}(s)ds\)\sqrt{\sN(w_r)\sD(w_r,\tilde w_r)}\nonumber\end{align}
for some
constant $\ell_2$. Then taking all possible $\gamma_{[r,t)}\in \cM[r,t)$ and $t\geq t'$, we have $$\kappa_{\th^*}(t')\leq \ell_2\(1+\int_{t'}^T\kappa_{\th^*}(s)/{\sqrt{T-s}}ds\)\text{ for all $t'\in[r,T]$}.$$
Note that  $\kappa_{\th^*}(T)\leq \ell_2\kappa_\zeta$, we have  $\kappa_{\th^*}(t)\leq \beta_2(T)(1+\kappa_\zeta)$ for some $\beta_2(T)$ i.e. $\th^*$ is P-Lipschitz with  Lipschitz constant $\beta_2(T)(1+\kappa_\zeta)$.

By  Lemma \ref{lemma-1} and \eqref{thtoTh},  $\Th^*(t,w_t,x)$ is horizontally and vertically differentiable for each fixed $x$ (thus P-continuous and  $\a$-differentiable) with
\begin{align}\label{verify}&\partial^\a\Th(t,w_t,x)\\
&=\int_{\dbR^n}\Psi_1(t,x;T,y)\partial^\a\hat h(t,w_t,y)dy+\int_{\dbR^n}\frac{d}{dt}[\Psi_1(t,x;T,y)]\hat h(t,w_t,y)dy\nonumber\\
&\q-\lan b(t, w_t,y), \theta(t, w_t,y)\ran-g(t,w_t,x)+\int_t^T\int_{\dbR^n}\Psi_1(t,x;s,y) \partial^\alpha\hat g(t,w_t,s,y)dy ds\nonumber\\
&\q+\int_t^T\int_{\dbR^n}\partial_t\Psi_1(t,x;s,y) \hat g(t,w_t,s,y)dy ds.\nonumber\\
&=-\lan b(t, w_t,y), \theta(t,w_t,y)\ran- g(t,w_t,x)-\cA\Th(t,w_t,x).\nonumber
\end{align}
In the last step, we used \eqref{lapaceTh}.
Moreover, when $t\rightarrow T^-$, it is easy to see that $\Th(T,w,x)=\lim_{t\rightarrow T^-}\Th(t,w_t,x)=h(w,x).$
Putting the above together, we conclude that $\Th$ is the classical solution of \eqref{PDPDE222}.

(5) {\it Some uniform estimates}. Finally, we prove the boundedness of $D_x\Th^*$ and  $D_x^2\Th^*$. By \eqref{repDTh}, for all $t\in[0,T]$, it follows that
\begin{align*}\sup_{(r,w_r,x)\in\dbM\times\dbR^n,r\geq t}\!\!\!\!\!\!\!\!|D_x\Th^*(r,w_r,x)|\leq \ell_3\(1+\int_t^T\sqrt{T-s}\!\!\!\!\!\!\!\!\sup_{(r,w_r,x)\in\dbM\times\dbR^n,r\geq s}\!\!\!\!\!|D_x\Th^*(r,w_r,x)|ds\)
\end{align*}
for some absolute constants $\ell_3$.
By Gronwall's inequality, there exists a constant $\beta_3(T)$ such that $|D_x\Th^*|\leq \beta_3(T)$.  By \eqref{Dypsi} and \eqref{thtoTh}, using integration by parts, we have \begin{align}\label{repD2Th}&D_{x_i,x_j}\Th^*(t,w_t,x)=\int_{\dbR^n}\Psi_1(t,x;s,y)\\
&\qq\q\times\[\int_{\cM}\[(D_{y_j}+\rho(t,x;s,y))(D_{y_i}+\rho(t,x;s,y))h(\tilde  w,y)\]\Psi_2(t,w_t;T,d\tilde w) \]dy\nonumber\\
&\q +\int_t^T \int_{\cM_s}\[\int_{\dbR^n}D_{x_i}\Psi_1(t,x;s,y)\Psi_2(t,w_t;s,d\tilde w_s) \nonumber \\
&\qq\qq\qq\times \[(D_{y_j}+\rho(t,x;s,y))(\lan b(s,\tilde w_s,y), \theta^*(s,\tilde w_s,y)\ran+g(s,\tilde w_s,y))\] dy ds.\nonumber\end{align}	
By \eqref{rhoesti},	it follows than  for any $t\in[0,T]$,
$$\sup_{(r,w_r,x)\in\dbM\times\dbR^n,r\geq t}\!\!\!\!\!\!\!\!\!\!\!|D_{x_ix_j}\!\Th^*(r,\!w_r,\!x)|\!\leq \!\ell_4\!\(1\!+\!\int_t^T\sqrt{T-s}\!\!\!\!\!\!\!\!\sup_{(r,w_r,x)\in\dbM\times\dbR^n,r\geq s}\!\!\!\!\!\!\!\!|D_{x_i,x_j}\Th^*(r,w_r,x)|ds\).$$ Then Gronwall's inequality yields that $D_x^2\Th^*$ is bounded uniformly with a constant $\beta_4(T)$. The proof is complete.
\end{proof}

We remark that the two constants $\beta_2(T)$ and $\beta_3(T)$ depend on $T$ and  are bounded if $T$ belongs to a finite horizon. Until now, we developed a theory for path-dependent PDEs. We are ready to apply the idea
to path-dependent equilibrium HJBs.

\subsection{Path-Dependent Equilibrium HJBs and
Local Optimality} 			
In this subsection, we focus on the  existence and uniqueness of the solution to  path-dependent equilibrium HJBs in \eqref{PDEHJB-00}. The following is the main result of this section.

\begin{theorem}\label{PDEHJBThm}
Suppose that $\zeta$ is P-continuous and  P-Lipschitz and that Assumption \ref{A-3} holds. The path-dependent system of equilibrium  HJBs in \eqref{PDEHJB-00} admits a  unique classical solution $\Th^\zeta$
such that $\partial^\a\Theta^\zeta(\t;t,w_t,x)$ and $D_{x}\Theta^\zeta(\t;t,w_t,x)$ are P-continuous for each $x$.
Moreover, $D_x\Th^\zeta(t,w_t,x)$ is P-Lipschitz with P-Lipschitz constant $\beta_2(T)(1+\kappa_\zeta)$ 
and
$| D_x\Theta^\zeta(t,w_t,x)|+| D_x^2\Theta^\zeta(t,w_t,x)|\leq  \beta_5(T)$
for some  constants $\beta_2(T),\beta_5(T)>0$ independent of $\zeta.$
\end{theorem}

\begin{proof} Similar to Proposition \ref{fundamental}, the proof is based
on a fixed point theory.
Given a $\gamma_{T-\d}\in\cM_{T-\d}$ and a $\theta_1\in \sC^{0,0}(\dbM_{\gamma_{T-\d}}\times\dbR^n)$, for each $\tau\in[0,T]$, we solve the following path-dependent  PDE \bel{PDEHJB-01}\left\{\ba{ll}\ns\ad(\partial^\a+\cA)\Theta(\tau; t,w_t,x)+\wt\cH(\t;t,w_t,x, D_x\Theta(\t;t,w_t,x),\theta_1(t,w_t,x))=0,\\ [2mm]
\ns\ad \Theta(\t;T,w,y)=h(\t; w,y).\ea\right.\eel	
By Proposition \ref{defvarphi0}, all the conditions in Proposition \ref{fundamental} are fulfilled for each $\tau\in[0,T]$. Then for  $\tau\in[0,T]$, there exists a unique solution $\Th_1(\tau;\cdot)$ with
$$\ba{ll}\ad\Th_1(\tau;t,w_t,x)=\int_{\cM\times \dbR^n}h(\t;\tilde w, y)\Psi(t,w_t,x;T,d\tilde w,y)dy\\
\ns\ad\q+\!\int_t^T\!\!\!\int_{\cM_s\times \dbR^n}\!\!\!\!\!\!\wt \cH(\t;s,\tilde w_s,x, D_x\Theta_1(\t;s,\tilde w_s,x),\theta_1(s,\tilde w_s,x))\Psi(t,w_t,x;s,d\tilde w_s,y)dyds.\ea$$
Let $ \th_2(t,w_t,x) :=D_x\Th_1(t;t,w_t,x)$. Repeating such process, we get a sequence $\{(\th_i,\Th_i)\}$ such that $ \th_{i+1}(t,w_t,x) :=D_x\Th_i(t;t,w_t,x)$ for $i=1,2,\ldots$. Similar to the proof of Proposition \ref{fundamental}, one can prove there exists a limit $(\th^*,\Th^*)$ satisfying $ \th^*(t,w_t,x) =D_x\Th^*(t;t,w_t,x)$ and  $D_x\Th^*(t;t,w_t,x)$, $D^2_x\Th^*(t;t,w_t,x)$, $\partial^\a \Th^*(t;t,\!w_t,x)$ are P-continuous for each $x$ and is continuous  w.r.t. $x$ for each $(t,w_t)\in\dbM_t$. Moreover, $D_x\Th(t;t,w_t,x)$ is P-Lipschitz with P-Lipschitz constant $\beta_2(T)(1+\kappa_\zeta)$.
\end{proof}

The solution $\Th^*(\tau;t,w_t,x)$ is the cost function when $u(\cdot)$ is taken as the equilibrium strategy.
Therefore, we can define the second mapping for HJB step by
letting $\cT_2:\Upsilon_{\mu_0}\mapsto\sU$ by
$\cT_2(\zeta)=u^\zeta$ where $u^\zeta(t,w_t,x)=\varphi\big(t,w_t(t^-),x,D_x\Th^\zeta(t;t,w_t,x)\big)$ for all $(t,w_t,x)\in\dbM\times\dbR^n$.
Different from \cite{Yong2012} and \cite{Mei2020},  our feedback strategy is path-dependent. 
Similar to \cite{Mei2020}, let us show that the feedback strategy satisfies a version of local optimality in the following proposition.

\begin{proposition}\label{localoptimalityustar} Given $\zeta$, let $\Th$ be the solution to the path-dependent equilibrium  HJBs in \eqref{PDEHJB-00}.
For any $u_0\in U$, it follows that
\bel{localopti} \begin{array}{ll}
\displaystyle \limsup_{\e\rightarrow 0^+}\frac1\e \[J(t;t,w_t,x;u^*)-J(t;t,w_t,x;u_0|_{[t,t+\e)}\oplus u^*|_{[t+\e,T]})\]\leq 0,\\
\hbox{\rm
where }
u_0|_{[t,t+\e)}\oplus u^*|_{[t+\e,T]}(t,w_t,x):=\left\{\ba{ll}\ns\ad  u_0,\qq\qq \text{ for } s\in [t,t+\e);\\
\ns\ad u^*(t,w_t,x),\qq \text{ for } s\in [t+\e,T].\ea\right. \end{array} \eel
\end{proposition}	

The proof is based on a functional It\^o's formula for $\Th(\tau;t,\a_t, X(t))$ in the following lemma.

\begin{lemma}\label{itofunctional} Given $\zeta,$ let $X$ be the solution to the following SDE,
\bel{SDE-200} dX(t)=\Bb(t,\alpha(t^-),X(t),\zeta(t,\a_t);u_0)dt+\BBsi(t, X(t))dB(t).\eel
Then
\begin{align}\label{hito}
&\dbE\Big[\Th(\tau;t+\e,\a_{t+\e},X(t+\e))\Big|\a_t=w_t,X(t)=x\Big]-\Th(\tau;t,w_t,x)\\
&=\e(\partial^\a+\cA)\Th(\tau;t,w_t,x)+\e\lan \Bb(t,w_t(t^-),x,\zeta(t,w_t),u_0),D_x\Th(\tau;t,w_t,x)\ran+o(\e).\nonumber
\end{align}
\end{lemma}

\begin{proof} Let $N_\e$ be the number of transitions of $\a(\cdot)$ in $[t,t+\e)$. Define the following events $A_\e(0)=I(N_\e=0),~A_\e(1)=I(N_\e=1),~A_\e(2)=I(N_\e\geq 2).$
One have $\dbP(A_\e(1))=\lambda\e+O(\e^2)\text{ and }\dbP(A_\e(2))=O(\e^2).$
Taylor's expansion yields that  for some $\xi$ on the line segment connecting $x$ and $X(t+\e)$,
\begin{align}\label{approITO}&\dbE_t\Th(\tau;t+\e,\a_{t+\e},X(t+\e))\!-\!\Th(\tau;t,w_t,x)
\!=\!\dbE_t\[\Th(\tau;t+\e,\a_{t+\e},x)-\Th(\tau;t,w_t,x)\]\\
&\q+\dbE_t\[D_x\Th(\tau;t+\e,\a_{t+\e},x)(X(t+\e)-x)\]\nonumber\\
&\q+\frac12\dbE_t\[\lan X(t+\e)-x,D_{xx}^2\Th(\tau;t+\e,\a_{t+\e},\xi)(X(t+\e)-x)\ran\]\nonumber\\
&=\dbE_t\[\(\Th(\tau;t+\e,w_{t}^\e,x)-\Th(\tau;t,w_t,x)\)I(A_\e(0))\]\nonumber\\
&\q+\dbE_t\[\(\Th(\tau;t+\e,\a_{t+\e},x)-\Th(\tau;t,w_t,x)\)I(A_\e(1))\]+O(\e^2)\nonumber\\
&\q+\e\[\dbE_tD_x\Th(\tau;t+\e,\a_{t+\e},x)\]\Bb(t,w_t(t^-),x,\zeta(t,w_t),u_0)\nonumber\\
&\q+\dbE_t\[D_x\Th(\tau;t+\e,\a_{t+\e},x)\int^{t+\e}_t \Bb(s,w_s(s^-),X(s),\zeta(s,
\a_s),u_0)\nonumber\\
&\qq\qq-\Bb(t,w_t(t^-),x,\zeta(t,w_t),u_0)ds\]\nonumber\\
&\q +\dbE_t\[\(D_x\Th(\tau;t+\e,\a_{t+\e},x)-\partial_x\Th(\tau;t,w_{t},x)\)\int_t^{t+\e}\sigma(s,x)dW(s)\]\nonumber\\
&\q+\frac12\dbE_t\[\lan X(t+\e)-x,D_{xx}^2\Th(\tau;t+\e,\a_{t+\e},\xi)(X(t+\e)-x)\ran\]\nonumber\\
&=\e\[\partial^H\Th(\tau;t,w_t,x)+\sum_{j=1}^m\partial_j^V\Th(\tau;t,w_t,x)q_{ij}\nonumber\\
&\qq+\lan D_x\Th(\tau;t,w_t,x),\Bb(t,w_t(t^-),x,\zeta(t,w_t),u_0)\ran+\cA\Th(\tau;t,w_{t},x)\]+o(\e).\nonumber
\end{align}
The last equality holds because
$\e^{-1/2}(X_{t+\e}-X_t)$ converges weakly to a normal distribution with mean $0$ and variance $\sigma^\top(t,x)\sigma(t,x)$, $\sup_\e\dbE|\e^{-1/2}(X_{t+\e}-X_t)|^4\leq L$ and $D_{xx}^2\Th(\tau;t+\e,\a_{t+\e},\xi)$ is bounded and converges $D_{xx}^2\Th(\tau;t,w_t,\xi)$ almost surely. Therefore we have
\begin{align*}&\lim_{\e\rightarrow 0^+}\frac1\e\dbE_t\[\lan X(t+\e)-x,D_{xx}^2\Th(\tau;t+\e,\a_{t+\e},\xi)(X(t+\e)-x)\ran\]\\
&\q=(\partial^\a+\cA)\Th(\tau;t,w_t,x)+ \lan \Bb(t,w_t(t^-),x,\zeta(t,w_t),u_0),D_x\Th(\tau;t,w_t,x)\ran.\end{align*}
The proof is complete.
\end{proof}

\begin{proof}[Proof of Proposition \ref{localoptimalityustar}] 	Let $ X$ be the solution under strategy $u^*$ and  $\bar X$ be the solution under strategy $u_0|_{[t,t+\e)}\oplus u^*(\cdot)|_{[t+\e,T]}$. Note that $$\Th(\t;t,w_t,x)=J(\t;t,w_t,x;u^*(\cdot)),$$ we have
$$\ba{ll}\ds J(t;t,w_t,x;u_0|_{[t,t+\e)}\oplus u^*|_{[t+\e,T]})- J(t;t,w_t,x;u^*(\cdot))\\
\ns\ds=\dbE_{t,w_t,x}\[\int_t^{t+\e}g(t;s,\a_s,\bar X(s), u_0)ds+\Th(t;t+\e,\a_{t+\e},\bar X(t+\e))\]\\
\ns\ds\q-\dbE_{t,w_t,x}\[\int_t^{t+\e}g(t;s,\a_s,X(s), u^*(s,\a_s,X(s)))ds+\Th(t;t+\e,\a_{t+\e},X(t+\e))\].\ea$$
Next, by the regularity for $\Th(\tau;t,w_t,x)$ in Theorem \ref{PDEHJBThm} and by virtue of estimates obtained in Lemma \ref{itofunctional},
$$\ba{ll}\ns\ad  \frac1\e\[J(t;t,w_t,x;u_0|_{[t,t+\e)}\oplus u^*|_{[t+\e,T]}(\cdot))-J(t;t,w_t,x;u^*(\cdot))\]\\
\ns\ad= \frac1\e\dbE\[\int_t^{t+\e}g(t;s,\a_s,\bar X(s), u_0)-g(t;s,\a_s,X(s), u^*(s,\a_s,X(s)))ds\]\\
\ns\ad\q+ \frac1\e\[\dbE\[\Th(t;t+\e,\a_{t+\e},\bar X(t+\e))-\Th(t;t+\e,\a_{t+\e},X(t+\e))\]\\
\ns\ad =g(t;t,w_t,x, u_0)- g(t;t,w_t,x, u^*(t,w_t,x))\\
\ns\ad\q+\lan D_x\Th(t;t,w_{t},x),b(t,x,w_t;u_0)-b(t,x,w_t;u^*(t,w_t,x))\ran+o(1)\geq o(1).\ea $$
In the last step, we have used the definition of $$u^*(t,w_t,x)=\psi(t,w_t(t^-),x,D_x\Th(t;t,w_t,x))$$ defined in \eqref{defpsi}. This verifies the
local optimality in \eqref{localopti}.
\end{proof}

We remark that the local optimality \eqref{localopti} yields the optimality of $u^*$ by dynamic programming principle if the cost functional is exponential discounting. Therefore, our result also solves the time-consistent optimal control problem with path-dependence on a switching environment.
As the HJB-step is completed,  we are ready to combine the HJB-step and the FP-step  for our
MFG problem.

\section{MFGs with Switching  under General Discounting Costs}\label{sec: mfg}

In this section, we
solve our main problem as stated in \eqref{cost} corresponding to SDE \eqref{SDE-1}. Our ultimate goal is to find
a closed-loop equilibrium  $u\in \sU$  as previously mentioned. Recalling  $\cT^{\mu_0}_1$ and $\cT_2$ from Section \ref{sec:fp} and Section \ref{sec: HJB-step} respectively, let us proceed with the  definition of a closed-loop equilibrium.

\begin{definition}\label{defofequilibrium}
\rm Given a $\mu_0\in\sP_2(\dbR^n)$,
a feedback strategy	
$u\in\sU$ is called a {\it closed-loop equilibrium} for initial $\mu_0$ if
$u=\cT_2\circ\cT^{\mu_0}_1(u).$
\end{definition}

  Our goal is  to obtain
  the existence and uniqueness of an equilibrium strategy.  It suffices to show $\cT_2\circ\cT^{\mu_0}_1$ is a contraction.  To achieve this,  we proceed with some {\it a priori} estimates.

\begin{lemma}\label{prioriestimateTh}
Suppose Assumption \ref{A-3} holds and $\zeta\in \Upsilon_{\mu_0}$ is P-continuous and P-Lipschitz.
Let $\Theta_i$ be the solution of \eqref{PDEHJB-00} using $\zeta_i$ for $i=1,2$. Then there exists a  constant $ \beta_5(T)>0$ such that
$$| D_x\Theta_1(\tau;t,w_t,x)-D_x\Theta_2(\tau;t,w_t,x)|\leq  \beta_5(T)d_\Upsilon(\zeta_1,\zeta_2).$$
Especially, $\beta_5(T)$ is bounded if $T$ is bounded.
\end{lemma}

\begin{proof}
 Write $$c_1=\sup_{\zeta_1,\zeta_2,\tau,(w,x)\in\cM\times\dbR^n}\frac{|h^{\zeta_1}(\t,w,x)-h^{\zeta_2}(\t,w,x)|}{d_\Upsilon(\zeta_1,\zeta_2)}.$$  Note that given a $\d>0$,  
 by \eqref{repDTh}, we have
\begin{align*}&
\ds|D_x\Theta_1(\t;t,w_t,x)-D_x\Theta_2(\t;t,w_t,x)|\\
&\leq \int_{\cM\times\dbR^n}|(D_y+\rho(t,x;s,y))(\Bh^{\zeta_1}(\t;\tilde w,y)-\Bh^{\zeta_2}(\t;\tilde w,y))|\Psi(t,w_t,x;T,d\tilde w,y)dy\\
&\ds\q+\int_t^T\int_{\cM_s\times\dbR^n}|\cH^{\zeta_1}(\t;s,\tilde w_s, y, D_y\Theta_1(\t;s,\tilde w_s,y),D_y\Theta_1(s;s,\tilde w_s,y))\\
&\ds\q-\cH^{\zeta_1}(\t;s,\tilde w_s, y, D_y\Theta_2(\t;s,\tilde w_s,y),D_y\Theta_1(s;s,\tilde w_s,y))| \\
&\ds\qq\times\frac{|y-x|}{s-t}\Psi(t,w_t,x;s,d\tilde w_s,y)dyds\\
&\ds\q+\int_t^T\int_{\cM_s\times\dbR^n}|\cH^{\zeta_1}(\t;s,\tilde w_s, y, D_y\Theta_2(\t;s,\tilde w_s,y);D_y\Theta_1(s;s,\tilde w_s,y))\\
&\ds\q -\cH^{\zeta_2}(\t;s,\tilde w_s, y, D_y\Theta_2(\t;s,\tilde w_s,y),D_y\Theta_2(s;s,\tilde w_s,y))| \\
&\ds\qq\times\frac{|y-x|}{s-t}\Psi(t,w_t,x;s,d\tilde w_s,y)dyds\\
&\ds\leq L\Big[d_\Upsilon(\zeta_1,\zeta_2)(1\!+\!c_1)+\!\int_t^T\!(s-t)^{-\frac12}\!\int_{\cM_s\times\dbR^n}\!\!\!\(|D_y\Theta_1(\t;s,\tilde w_s,y)\!-\!D_y\Theta_2(\t;s,\tilde w_s,y)|\\
&\q\!+\!|D_y\Theta_1(s;s,\tilde w_s,y)\!-\!D_y\Theta_2(s;s,\tilde w_s,y)|\!+ \! \frac{d_\Upsilon(\zeta_1,\zeta_2)|y-x|}{(s-t)^\frac12}\Psi(t,w_t,x;s,d\tilde w_s,y)dyds\Big]\\
&\ds\leq L\Big[d_\Upsilon(\zeta_1,\zeta_2)(1+c_1)+ (T-t)^\frac12\!\!\!\!\!\!\!\!\!\!\!\!\!\! \sup_{\tau,r\geq T-\d, (r,\tilde w_r,y)\in\dbM\times \dbR^n}\!\!\!\!\!\!\!\!\!\!\!\!\!\! |D_y\Theta_1(\t;r,\tilde w_r,y)-D_y\Theta_2(\t;r,\tilde w_r,y)|\Big].
\end{align*}
Then it follows that
$$\ba{ll}\ad\sup_{\tau,t\geq T-\d,(t,w_t,x)\in\dbM\times\dbR^n} |D_y\Theta_1(\t;t,\tilde w_t,y)-D_y\Theta_2(\t;t,\tilde w_t,y)|\\
\ns\ad\leq C_3(T) \Big[d_\Upsilon(\zeta_1,\zeta_2)(1+c_1) +\sqrt{\d} \!\!\!\!\!\!\!\!\!\!\!\!\!\! \sup_{\tau,t\geq T-\d,(t,w_t,x)\in\dbM\times\dbR^n}\!\!\!\!\!\!\!\!\!\!\!\!\!\!\! \!\!\! |D_y\Theta_1(\t;t,\tilde w_t,y)-D_y\Theta_2(\t;t,\tilde w_t,y)|\Big] .\ea$$
Here $C_3(T)$ is a constant dependent  of $T$ but independent of $c_1,\zeta_i$. Take $\d=C_3^{-2}(T)/2$, then
\begin{align}\label{t-dtoT}&\sup_{\tau,t\geq T-\d,(t,w_t,x)\in\dbM\times\dbR^n}| D_x\Theta_1(\t;t,w_t,x)\!-\!D_x\Theta_2(\t;t,w_t,x)|\\
&\q\leq C_3(T)(1+c_1)(1-C_3(T)\sqrt{\d})^{-1}d_\Upsilon(\zeta_1,\zeta_2). \nonumber\end{align}
Because $\d$  is independent of $c_1$ and $\zeta_i$,  we can repeat
the above procedure on $[T-2\d, T-\d]$ using $\Th_i(T-\d, w_{T-\d},x)$ as a terminal condition. Especially, on the terminal, we have $$|\Th_1(\tau;T-\d,w_{t-\d},x)-\Th_2(\tau;T-\d,w_{t-\d},x)|\leq  C_3(T)(1+c_1)(1-C_3(T)\sqrt{\d})^{-1}d_\Upsilon(\zeta_1,\zeta_2).$$ Write $c_2= C_3(T)(1+c_1)(1-C_3(T)\sqrt{\d})^{-1}$. Then we have a similar result to \eqref{t-dtoT} by replacing $c_1$  with $c_2,$ i.e.
\begin{align*}&\sup_{\tau,t\geq T-2\d,(t,w_t,x)\in\dbM\times\dbR^n}| D_x\Theta_1(\t;t,w_t,x)\!-\!D_x\Theta_2(\t;t,w_t,x)|\\
&\q\leq C_3(T)(1+c_2)(1-C_3(T)\sqrt{\d})^{-1}d_\Upsilon(\zeta_1,\zeta_2). \nonumber\end{align*}
Because the choice of $\d$ is independent of $c_i$, we
repeat the procedure until $t=0$ in finite induction. Then there exists a
$\beta_5(T)>0$
such that
$| D_x\Theta_1-D_x\Theta_2|\leq \beta_5(T)d_\Upsilon(\zeta_1,\zeta_2). $
The proof is complete.
\end{proof}

Now, we are ready to present the main result of this paper.

\begin{theorem}\label{mainthm2221}
Suppose Assumption \ref{A-3} holds.
Given $\mu_0\in\sP_2(\dbR^n)$, if \begin{equation}\label{smalnness} L_\psi \beta_5(T)\sqrt{\beta_1(T)}, ~\beta_2(T)\sqrt{\beta_0(T)}<1,\end{equation} then $\cT_2\circ\cT_1^{\mu_0}$
admits a unique fixed point in $\sU$, which  is the closed-loop equilibrium in the form of
$$u^*(t,w_t,x)=\psi(t,w_t(t^-),x,D_x\Theta^*(t;t,w_t,x)),$$
where $\Theta^*(\t;t,w_t,x),$ is the solution to \eqref{PDEHJB-00}. 
\end{theorem}

\begin{proof}
For $i=0,1,\ldots$, let $\zeta_{i+1}=\cT_1^{\mu_0}(u_i)$ and $u_{i+1}=\cT_2(\zeta_{i+1})$.  Let $\Theta_i$ be the  solution to \eqref{PDEHJB-00} w.r.t. $\zeta_i$. Then
$u_i(t,w_t,x)=\psi(t,w_t(t^-),x,D_x\Theta_i(t;t,w_t,x)).$
By Theorem \ref{PDEHJBThm}, we have
\bel{uniformDDu}|D_xu_i(t,w_t,x)|\leq  L(1+|D_x^2\Th_i(t;t,w_t,x)|)\leq L\big(1+  C_2(T)\big).\eel
By \eqref{zetauniform}, it follows that
$$\ba{ll}\ad|u_{i+1}(t,w_t,x)-u_i(t,w_t,x)|\leq L_{\psi}|D_x\Th_{i+1}(t;t,w_t,x)-D_x\Th_i(t;t,w_t,x)|\\
\ns\ad\leq L_\psi \beta_5(T)d_\Upsilon(\zeta_1,\zeta_2)\leq  TL_\psi \beta_5(T)\sqrt{\beta_1(T)}\!\!\!\!\sup_{(t,w_t,x)\in\dbM \times\dbR^n}\!\!\!\!{|u_{i}(t,w_t,x)-u_{i-1}(t,w_t,x)|}.\ea$$
Our assumption yields that the mapping $\cT_2\circ\cT_1^{\mu_0}$ is a contraction which admits a fixed point    $u^*\in
\sC^{0,0}(\dbM\times\dbR^n:\dbR^n)$. By \eqref{zetauniform} and  Proposition \ref{fundamental}, for the P-Lipschitz constant $\kappa_{u_i}$ of $u_i$, we have
$$\kappa_{u_{i+1}}\leq \beta_2(T)(1+\kappa_{\zeta_{i+1}})\leq  \beta_2(T)\big(1+(1+ \kappa_{u_i})\sqrt{\beta_0(T)}\big).$$
Given $\beta_2(T)\sqrt{\beta_0(T)}<1$, $\sup_i\kappa_{u_i}<\infty $ and therefore $u^*$ is P-Lipschitz.
Moreover, by the uniform estimate in \eqref{uniformDDu}, we have $u^*\in\sU$, which is the closed-loop equilibrium. The proof is complete.
\end{proof}

Now we have derived the existence and uniqueness of a closed-equilibrium for our MFG problems in a switching environment under general discounting costs. In our main theorem, the smallness assumption \eqref{smalnness}  is assumed such that the map is a contraction.   As we mentioned in Remark \ref{zetauniform} and Lemma \ref{prioriestimateTh}, $\beta_0(T)$ and $\beta_1(T)$ are small when $T$ is small, and $\beta_5(T)$ and $\beta_2(T)$ are bounded when $T$ is bounded. Therefore  when $T$ is small, \eqref{smalnness} holds  and our results follow.

Finally, let us return to the motivational example given in the introduction.   Theorem \ref{mainthm2221} and Proposition \ref{localoptimalityustar} lead to the following corollary directly.

\begin{corollary} For the motivational Example \ref{m-exm}, suppose all the assumptions in Theorem \ref{mainthm2221} holds. Then there exists a closed-loop equilibrium $u^*\in\sU$ which is the fixed-point of the HJB-FP method such that the following local optimality holds,
\begin{equation}\label{optmaNash2}\limsup_{\e\rightarrow 0^+}\frac{1}{\e}\(\wt J(t;t,x;\eta^{\xi,u^*},u_0|_{[t,t+\e)}\oplus u^*|_{[t+\e,T]}) -\wt J(t;t,x;\eta^{\xi,u^*},u^*)\)\leq 0\end{equation}
for any  $u_0\in U$.
\end{corollary}

Note that in the above, we used ``$\leq$"  because
of the maximization of
the profit functional in the motivational example.   Different from \eqref{optmaNash},   \eqref{optmaNash2} only requires that  the closed-loop equilibrium  to be  asymptotically optimal at a small time horizon $[t,t+\e)$  if all the agents follow the pre-committed closed-loop equilibrium  afterward and the agent adopts the discounting factor $\tau$ as the decision time $t$. In fact,  the new closed-loop equilibrium  is also a generalization of the previous one, because   \eqref{optmaNash2} is sufficient for \eqref{optmaNash} by the dynamical programming principle if $\wt J$ in \eqref{costEx2} is independent of $\tau$.
As can be seen, our result
has extended the classical MFG theory  to conditional MFGs with  general discounting costs.

\section{Concluding Remarks}\label{sec:cr}
In this paper, we focused on a class of
MFGs in randomly switching environments
with general discounting costs. The problem has two main features, namely, path-dependence due to conditional mean-field interactions, and time-inconsistency due to general discounting costs. We studied  the calculus on the path space of the Markov chain and developed  a theory on path-dependent equilibrium HJB equations to overcome the  difficulties brought by the two new features. Finally,  we proved the existence of a unique closed-loop equilibrium  for the problem of interest and verified that the closed-loop equilibrium satisfies a local Nash's optimality instead.  As a result, this paper presents a substantial extension to classical MFG theory as well as the theory of  equilibrium HJB equations in a path-dependent setting (c.f., for example,  \cite{Yong2012,Mei2020} and references therein for previous works in the literature).
It is conceivable that
techniques developed in this work can also be used to treat
mean-field optimal control of switching diffusions.
For future works, a number of problems may be considered further.
One may wish to consider the Hamiltonian that is not in the separable
form; one may also consider the diffusion coefficient depending on the control process. All of these deserve careful thought and further consideration.

\bibliographystyle{siamplain}

\end{document}